# Theoretical and computational investigations of superposed interacting affine and more complex processes

Hidekazu Yoshioka[1,*]

[1] Japan Advanced Institute of Science and Technology, 1-1 Asahidai, Nomi 923-1292, Japan
[*]Corresponding author: yoshih@jaist.ac.jp, ORCID: 0000-0002-5293-3246

**Abstract**
We theoretically and computationally investigate long-memory processes based on the Markovian lifts of affine jump-diffusion processes. A nominal superposition process consisting of an infinite number of interacting affine processes is considered, along with its finite-dimensional version and associated generalized Riccati equations. We propose a splitting scheme suited to the Markovian lifts where jump and diffusion parts are dealt with separately based on recently developed exact discretization methods. We examine the computational performance of the scheme through comparisons with the analytical results. We also numerically investigate a more complex model arising in the environmental sciences and some extended cases in which superposed processes belong to a class of nonlinear processes that generalize affine processes.

***Statements & Declarations***
**Fundings** This study was supported by the Japan Society for the Promotion of Science (KAKENHI No. 22K14441 and 25K00240) and Japan Science and Technology Agency (PRESTO No. JPMJPR24KE). We thank the anonymous reviewers for their helpful comments.
**Competing interests** The authors have no relevant financial or nonfinancial interests to disclose.
**Data availability** The data will be made available upon reasonable request from the corresponding author.
**Acknowledgments** N.A.
**Declaration of generative AI in scientific writing** The authors did not use generative AI for the scientific writing of this manuscript.

## 1. Introduction

### 1.1 Research background

Stochastic processes are mathematical models that describe each time-dependent phenomenon by focusing on its fluctuations or driving noise, enabling us to study their dynamics theoretically and computationally [1-3]. In particular, a stochastic process with an autocorrelation function that asymptotically decays at an algebraic speed ($h^{-a}$ with some $a > 0$ for time lag $h >> 1$) is referred to as a long-memory process [4]. Long memory processes arise in various problems such as virtual currency fluctuations [5], climate and precipitation dynamics [6,7], natural language sequences [8], and dynamics of the solar coronal magnetic field [9]. The appearance of a long memory in a stochastic process is often regarded to be non-Markovian, indicating a certain path dependence and persistent correlation in the process [e.g., 10-13]. Directly dealing with a non-Markov process is a challenging task because of the lack of Markovian nature, just as the name suggests; indeed, the modern stochastic calculus largely depends on Markov processes [e.g., 14,15].

A breakthrough to resolve this difficulty, called Markovian lifts, was investigated by Carmona and Coutin [16], focusing on fractional Brownian motion, whose main concept was to rewrite a non-Markov process to (possibly a continuum of infinitely) many Markov processes, to which common mathematical tools in stochastic calculus such as the classical Itô's formula apply. Each Markov process in the Markovian lifts has a distinctive reversion speed, i.e., timescale, covering a broad range of timescales involved in the original non-Markov process as a whole, and the coexistence of multiple timescales reproduces long memory. That is, aggregating Markov processes with distinct timescales can lead to a non-Markov process. Markovian lifts have been effectively applied to various stochastic process models with memory, such as the autoregressive processes [17,18], fractional Ornstein–Uhlenbeck processes [19-21], Volterra processes arising in economics [22,23] and physics [24,25], and superposition of continuous-state branching processes and their variants [26,27]. The Erlang polynomial approach [28], which decomposes an Erlang memory into a sequence of exponential memories, can be considered a discrete version of Markovian lifts. An application to a mean-field game [29] was also conducted, which was a deterministic model that used a common principle to superpose processes with multiple timescales to obtain persistent memory.

Despite the success of Markovian lifts in theory, their numerical computation is not always easy because of the increase in problem dimensions from one to infinity. Approximating the Markovian lifts is often used in applications so that only a finite number of Markov processes appear; this type of truncation approach has been used for Volterra stochastic differential equations (SDEs) that are formulated as stochastic convolutions between non-exponential kernel functions and Lévy(-driven) processes [30-33]. By contrast, computational investigations of superposition processes, which are generated by some spatial integration (i.e., superposition) of measure-valued processes, driven by Lévy bases (i.e., infinite-dimensional Lévy processes) have been addressed for limited cases, such as the finite difference method for the superposition of Ornstein–Uhlenbeck processes [34], kernel-based methods for trawl processes [35,36], the finite difference method for ambit fields [37], and the Fourier approximation method for ambit fields [38]. Curato et al. [39] studied spatiotemporal ambit fields and their applications to the learning problem of raster data defined on grids.

Recently, the authors considered a superposition process based on square-root diffusion, focusing on water quality dynamics in streamflow environments [40]. In particular, computational studies on superposition processes based on interacting (i.e., mutually correlated) affine processes are still rare, except for that of Yoshioka [27], who discussed a fish migration phenomenon; however, the processes studied were purely jump-driven processes, and there still exist many other superposition processes that could potentially arise in applications. Moreover, nonlinear superposition processes beyond linear and affine Markov processes have not been computationally examined to the best of the author's knowledge but have been suggested to be possible, at least for certain Volterra processes [41]. These research backgrounds motivated this study explained below.

### 1.2 Aim and contribution

This study has two aims **A-B**. A greater weight is placed on aim **A** because aim **B** stems from it.

**A.** Analysis and computation of a class of superposition processes based on interacting affine processes.

**B.** Numerical investigations of more complicated cases, including nonlinear ones.

The contributions to achieving each aim are as follows.

Regarding aim **A**, we consider a nominal model as a non-Markov process arising from the superposition of a continuum of infinitely many affine jump-diffusion processes, and hence is considered infinite-dimensional, which interact with each other. The nominal model is an infinite-dimensional version of existing finite-dimensional affine processes such as multivariate Hawkes [42-46] and Hawkes jump-diffusion processes [47-49]. Infinite-dimensional versions of self-exciting pure jump processes have been considered by Agathe-Nerine [50], focusing on neural activities modeled through dense random graphs, which can be considered as a version of the nominal model with specific connectivity among the processes to be superposed. A field-theoretic method has been implemented for an infinite-dimensional system arising from a nonlinear Hawkes process to recursively determine its cumulants [51], which is similar to ours because of the lifting of the dimension despite the qualitative difference between the driving noise processes.

The affine nature of the processes to be superposed in the nominal model not only allows for self-exciting phenomena owing to interactions but also derivations of cumulants analytically or through a generalized Riccati equation, the latter being a partial integro-differential equation. We study the convergence of a finite-dimensional system arising from the truncated Markov lifts of the nominal model. Moreover, we study the generalized Riccati equation because it appears as a unique partial integro-differential equation whose solution existence is a non-trivial issue. We demonstrate that the generalized Riccati equation admits a unique solution obtained by a variation of constants formula.

The infinite-dimensional nature of the generalized Riccati equations suggests a linkage between the proposed model and superprocesses, such as super-Brownian motions and related models [52,53, Chapter 4 in 54], which satisfy certain partial differential equations. The numerical computation of super-Brownian motion has not been well documented, except by Moro and Schurz [55], who developed a nonnegativity-preserving discretization method. Although the forms of the equations and their associated Riccati equations differed between our model and these super-processes, this study provides insights into both models because they are driven by space-time white noise processes.

The numerical computation of the nominal model is based on exact simulation methods for jump and square-root processes that are less biased than common numerical methods such as the Euler–Maruyama method because the latter has been proven to be non-convergent when the diffusion coefficient is large [56], which is a serious issue in applications. These numerical methods achieve consistency between conditional moment-generating functions and discretized processes [57,58]. We exploit the advantage that they can be applied to a full range of parameter values. Indeed, square-root diffusion, as the simplest affine process, has distinct behaviors for large and small diffusion constants such that the process touches the boundary in the former, while the latter does not (e.g., Chapters 1-2 in Alfonsi [59]). We propose an operator-splitting method that harmonizes with exact discretization methods, with which the nonnegativity of solutions to the nominal model is preserved. We also demonstrate that the exponential moments obtained by exact discretization methods agree well with those obtained from the generalized Riccati equation.

Regarding Aim **B**, we apply the proposed numerical method to the coupled system of two superposition processes representing water quantity-quality dynamics in a streamflow environment [40]. This model admits a closed-form average and covariance and has been applied to river discharge (jump-driven superposition process) and multiple water quality indices (diffusion-driven superposition process) [40]. This study is not intended to compare methods to deal with superposition processes but rather to examine a target methodology from both theoretical and computational standpoints. We apply exact discretization methods to this coupled model and demonstrate that they compute a model that is more complex than a nominal one. We also consider a constant elasticity of variance (CEV) model [60-62] that nonlinearly generalizes the square-root diffusion term in the nominal model; thus, the superposition approach is applied formally. We investigate the parameter dependence of this generalized model with respect to the regularity of diffusion coefficients. Through these test cases and application studies, we examine the applicability and limitations of the proposed numerical scheme based on exact discretization methods for long-memory processes. The outcomes of this study will contribute to both the modeling and computation of long-memory processes.

### 1.3 Structure of this paper

The remainder of this paper is organized as follows. **Section 2** presents and analyzes the nominal model used in this study. **Section 3** reviews the exact discretization methods used in this study and applies them to the proposed model. **Section 4** discusses the computational performance and applications of the numerical scheme. **Section 5** summarizes the study and presents perspectives. **Appendix** presents proofs of a lemma and propositions.

## 2. Nominal model

## 2.1 SDE formulation

We use a filtered probability space $\left(\Omega,\mathbb{F},(\mathbb{F}_t)_{t\in\mathbb{R}},\mathbb{P}\right)$. We consider a superposition process model that generates a continuous-time scalar process $X=(X_t)_{t\in\mathbb{R}}$. Its dynamics are assumed to be driven by two independent Lévy bases, which are pure jump and diffusion ones. Based on the superposition approach (e.g., Barndorff-Nielsen [63]), we consider that $X$ is formally an integration of a measure-valued process $x(\cdot)=\left(x_t(\cdot)\right)_{t\in\mathbb{R}}$ parameterized by the reversion speed $r\in(0,+\infty)$ generated by the probability measure $\rho$:

$$X_t=\int_0^{+\infty}x_t(\mathrm{d}r),\quad t\in\mathbb{R}. \tag{1}$$

For any Borel measurable set $S\subset(0,+\infty)$ and time $t\in\mathbb{R}$, the measure-valued process $x$ is assumed to follow the SDE

$$\underbrace{\mathrm{d}x_t(S)}_{\text{Increment}}=\int_{r\in S}\left(\underbrace{-rx_t(\mathrm{d}r)\mathrm{d}t}_{\text{Mean reversion}}+\underbrace{rb(r)\rho(\mathrm{d}r)\mathrm{d}t}_{\text{Drift}}+\underbrace{B(\mathrm{d}t,\mathrm{d}r)}_{\text{Diffusion}}+\underbrace{\int_0^{+\infty}zJ(\mathrm{d}t,\mathrm{d}r,\mathrm{d}z)}_{\text{Jump}}\right), \tag{2}$$

which can also be rewritten as follows (called nominal model in the rest of this paper):

$$\mathrm{d}x_t(\mathrm{d}r)=-rx_t(\mathrm{d}r)\mathrm{d}t+rb(r)\rho(\mathrm{d}r)\mathrm{d}t+B(\mathrm{d}t,\mathrm{d}r)+\int_0^{+\infty}zJ(\mathrm{d}t,\mathrm{d}r,\mathrm{d}z),\quad r>0. \tag{3}$$

Each term in (2) is explained as follows. The left-hand side is the increment of $x$ during a time increment $\mathrm{d}t$. The first term in the right-hand side represents the mean reversion of $x$ with the reversion speed $r$. The second term with a smooth function $b:(0,+\infty)\to[0,+\infty)$ is deterministic drift that is strictly bounded by a constant from above. The multiplication by $\rho(\mathrm{d}r)$ is to ensure that the right-hand side of (2) is measure-valued. The third term is Gaussian; a continuous increment with a Gaussian measure $B$ whose covariance satisfies

$$\mathbb{E}\left[B(\mathrm{d}t,\mathrm{d}r)B(\mathrm{d}t',\mathrm{d}r')\right]=\sigma^2r\delta(r-r')\delta(t-t')x_t(\mathrm{d}r) \tag{4}$$

using the noise intensity $\sigma\geq 0$ and the Dirac delta $\delta(\cdot)$. State-dependent covariance (4) resembles square root diffusion as the simplest affine diffusion process (Chapters 1–2 in Alfonsi [59]). Finally, the fourth term in (2) represents a jump; a discontinuous increment with the compensated version $\tilde{J}$ of $J$ is given as follows so that jumps contain external and self-exciting ones:

$$\tilde{J}(\mathrm{d}t,\mathrm{d}r,\mathrm{d}z)=J(\mathrm{d}t,\mathrm{d}r,\mathrm{d}z)-\left(\underset{\text{External factor}}{ra}+\underset{\text{Self-exciting factor}}{\mu_x(r)}\right)\rho(\mathrm{d}r)\nu(\mathrm{d}z)\mathrm{d}t. \tag{5}$$

Here, $\nu$ denotes a Lévy measure of a pure-jump process with positive jumps such that $M_1<+\infty$ ($M_k=\int_0^{+\infty}z^k\nu(\mathrm{d}z)$, $k=1,2,3,\ldots$). This condition is not restrictive in applications because it is satisfied by typical compound Poisson processes and pure jump subordinators (Lemma 2.14 in Kyprianou [64]). In (5), the intensity owing to the self-exciting jump is set as follows:

$$\mu_x(r) = \int_0^{+\infty} uA(r,u) x_t(\mathrm{d}u) \tag{6}$$

with $A:(0,+\infty)^2 \to [0,+\infty)$ and a constant $a \geq 0$. Unless otherwise specified, we assume that $A$ is bounded, nonnegative, and uniformly continuous in $(0,+\infty)^2$.

We compare technical aspects between the present and existing approaches for infinite-dimensional SDEs. Our nominal model is formulated based on the superposition approach [63] that has originally been discussed for a superposition of Ornstein–Uhlenbeck processes, and its solution $x_t(\mathrm{d}r)$ can be understood in a weak sense of an infinite-dimensional limit of a finite-dimensional one to be presented in **Section 2.3** (see also Iglói [26]). This is an intuitive approach where an integration of a measure-valued process is replaced by a summation of a system of usual SDEs in computational implementations. The superposition is conducted with respect to the parameter $r > 0$ whose inverse corresponds to timescale. This approach is different from those based on Hilbert spaces where a superposition is conducted with respect to a proper orthogonal basis and increasing the number of bases improves accuracy of the finite-dimensional system, which is a projection of the infinite-dimensional one [65,66]. In this view, our approach is closer to that for Volterra SDEs [30-33] where the discretization of a kernel corresponds to that for $r > 0$ in our model. The Volterra and superposition approaches are different from that based on Hilbert spaces because the former try to derive an infinite-dimensional process by some discretization in space, which is fundamentally different from the expansion by basis employed in the latter. The Volterra and superposition approaches are also technically different from each other because of superposing processes driven by common and independent noises, respectively. Finally, a common feature of the three approaches is that the analysis of the target process can be facilitated by studying the associated generalized Riccati equation.

Superposing measure-valued processes (2) as in (1) can yield a stochastic process that can capture various memories; however, the mathematical difficulty is that the superposition process $X$ is not always Markovian. In general, we have $\mathbb{E}[X_t | X_s] \neq \mathbb{E}[X_t | \mathbb{F}_s]$ ($t > s$) because the information obtained by collecting the entire $x_s(\cdot)$, i.e., $\mathbb{F}_s$, is richer than that obtained by knowing the value of $X_s$. From a quantitative viewpoint, the memory of $X$ is generated by the coexistence of processes having different timescales ranging from $r = 0$ to $+\infty$.

We set the following assumptions that are assumed in the rest of this paper:

***Assumption 1*** *There exists a constant* $\bar{A} > 0$ *such that*

$$\sup_{r,u>0} A(r,u) \leq \bar{A} < +\infty. \tag{7}$$

***Assumption 2***

$$M_1, M_2 < +\infty \quad \text{and} \quad \bar{A} M_1 < 1. \tag{8}$$

***Assumption 3*** *The probability measure* $\rho$ *admits a density. In addition, the density is positive and bounded at all* $r>0$ *and has an average.*

We consider that these assumptions are not restrictive in practice. In applications, particularly those arising in environmental sciences and related research fields, research interests are often long-term assessment where the target variable is based on some stationary processes. The self-exciting jumps are not highly frequent so that stationary average and variance exist (**Assumptions 1 and 2**), and the existence of density (**Assumption 3**) is not restrictive in applications because assuming $\rho$ as the probability measure of gamma distribution is able to generate a wide range of memory decays (see **Section 4**).

## 2.2 Related models

### 2.2.1 Specific cases

A few models represented by and related to nominal model (3) are presented. The simplest case covered by the nominal model is the superposition of Ornstein–Uhlenbeck processes ($b\equiv 0$, $\sigma=0$, $A\equiv 0$) [63] in finance and economics, and the superposition of diffusion processes with a linear differential generator ($a=0$, $A\equiv 0$) in network traffic (Chapter 4 in Iglói [26]). The two aforementioned models are superpositions of independent stochastic processes, whereas the superposition of continuous-time branching-type processes with mutual interactions ($\sigma=0$, $A\neq 0$) was proposed by Yoshioka [27] to study fish migration. Specifically, we can formally consider the case in which $A(r,u)=A_1(r,u)+A_2\delta(u-r)$ with $A_1$ is bounded and continuous and $A_2\geq 0$ is a constant. In this case, (6) becomes:

$$\mu_x(r)=A_2 r x_t(\mathrm{d}r)+\int_0^{+\infty} u A_1(r,u)x_t(\mathrm{d}u), \tag{9}$$

which decomposes the intensity of the self-exciting jumps into those owing to itself (first term) and aggregation (second term). Yoshioka [27] considered a superposition process without diffusion and $A_1, A_2$ as small positive constants.

The nominal model (3) assumes a time-independent coefficient $b$, while it can be stochastic and time-dependent and was recently proposed by Yoshioka and Yoshioka [40] to describe the concentration-discharge relationships in a streamflow environment. Their model was a couple of superposition processes in which the coefficient $b$ was replaced by a nonnegative long memory process. We consider that a greater extension of the nominal model would be possible where each coefficient, such as $b$, nonlinearly depends on the $X$ of some integration of $x(\mathrm{d}r)$, although such an extension does not allow for the application of the theory of the generalized Riccati equation. However, such a nonlinear model can be computed numerically by suitably discretizing each integral.

Lévy bases serve as generalizations of Gaussian white noise and Poisson random measures, which are considered infinite-dimensional versions of Lévy processes [67]. Models driven by generic Lévy bases can also be considered by replacing the noise $B$ and $J$ and/or adding state-independent terms driven by them.

### 2.2.2 Specific cases

We also consider a nonlinear generalization of the nominal model. The superposition approach assumes that the allowable nonlinearity with respect to the measure-valued process $x$ can be limited because each $x$ has a size of $\mathrm{d}r$ in both the average and variance of the superposition processes (e.g., Section 2 in Yoshioka [27]). For example, considering a CEV-type generalization [60-62] in the nominal model would be

$$\sigma^2 r\delta(r-r')\delta(t-t')x_t(\mathrm{d}r) \Rightarrow \sigma^2 r\delta(r-r')\delta(t-t')(X_t)^\gamma \tag{10}$$

with $\gamma \neq 1$. This case corresponds to an $X$-dependent coefficient $A$ (e.g., a superposed version of Qu et al. [68]) or multiplying a function of $X$ by a white-noise $B$. We arrive at the classical CEV model when there is no superposition, that is, when $\rho$ is a Dirac delta concentrated at a positive value. This nonlinear generalization is computationally discussed in **Section 4**.

## 2.3 Finite-dimensional version

A finite-dimensional version of the nominal model, in which the probability measure $\rho$ is replaced by an empirical measure, is proposed. This model is a building block for the numerical scheme proposed in this study, as well as an intuitive model for studying the generalized Riccati equation. The discretization in space used in this study is the quantile-based method presented by Yoshioka et al. [69]. In the remainder of this section, we assume that $\rho$ admits a density.

For a fixed degree of freedom $N \in \mathbb{N}$, the probability measure $\rho$ is approximated by the empirical one $\rho_N$:

$$\rho \approx \rho_N = \frac{1}{N}\sum_{i=1}^{N}\delta(r-r_i), \tag{11}$$

where $r_i$ is the $\frac{2i-1}{2N}$ th quantile level of $\rho$, i.e., $\int_0^{r_i}\rho(\mathrm{d}r) = \frac{2i-1}{2N}$ ($i = 1,2,3,...,N$). Similarly, set the $\frac{2i}{2N}$ th quantile level $\overline{r}_i$ of $\rho$ by $\int_0^{\overline{r}_i}\rho(\mathrm{d}r) = \frac{2i}{2N}$ ($i = 0,1,2,...,N$, $\overline{r}_N = +\infty$). We have $\overline{r}_{i-1} < r_i < \overline{r}_i$ and $\int_{\overline{r}_{i-1}}^{\overline{r}_i}\rho(\mathrm{d}r) = \frac{1}{N}$, $i = 1,2,3,...,N$, justifying the weight $\frac{1}{N}$ in (11). The accuracy of this quantile-based discretization is $\frac{1}{N}$ with respect to the supremum norm of the cumulative distribution (e.g., Proof of Proposition 2 in Yoshioka [70]).

We write $b_i = b(r_i)$ and $A_{i,j} = A(r_i, r_j)$. Substituting $S = (\overline{r}_{i-1}, \overline{r}_i)$ into (2) yields the system of Itô's SDEs that governs $\left\{x_t^{(i)}\right\}_{i=1,2,3,...,N} \approx \left\{x_t\left(\left(\bar{r}_{i-1},\bar{r}_i\right)\right)\right\}_{i=1,2,3,...,N}$ :

$$\begin{array}{l} \mathrm{d}x_t^{(i)} = -r_i x_t^{(i)}\mathrm{d}t + \dfrac{1}{N} r_i b_i \mathrm{d}t + \sigma\sqrt{r_i x_t^{(i)}}\mathrm{d}W_t^{(i)} \\ \quad + \displaystyle\int_{u=0}^{u=\frac{1}{N}(r_i a+\mu_i)} \int_{z=0}^{z=+\infty} zJ^{(i)}\left(\mathrm{d}t,\mathrm{d}u,\mathrm{d}z\right) \end{array}, \quad t \in \mathbb{R}, \quad i = 1,2,3,...,N \tag{12}$$

with

$$\mu_i = \sum_{j=1}^{N} r_j A_{i,j} x_t^{(j)} , \tag{13}$$

where $\left\{W^{(i)}\right\}_{i=1,2,3,...N}$ is a sequence of one-dimensional Brownian motions that are mutually independent, and $\left\{J^{(i)}\right\}_{i=1,2,3,...N}$ is a sequence of Poisson random measures with compensated measures that are mutually independent with the common compensation measure $\mathrm{d}u\nu\left(\mathrm{d}z\right)\mathrm{d}t$. We used the discretization $\int_{r\in\left(\bar{r}_{i-1},\bar{r}_i\right)}\int_0^{+\infty} zJ\left(\mathrm{d}t,\mathrm{d}r,\mathrm{d}z\right) \approx \int_{u=0}^{u=\frac{1}{N}(r_i a+\mu_i)} \int_{z=0}^{z=+\infty} zJ^{(i)}\left(\mathrm{d}t,\mathrm{d}u,\mathrm{d}z\right)$ for the last term in (12). The two noise sequences are assumed to be mutually independent. The discretized superposition process $X$ is then set to

$$X_t^{(N)} = \sum_{i=1}^{N} x_t^{(i)} . \tag{14}$$

### 2.4 On convergence and stationary solution

We now present remarks on the unique existence of a strong (pathwise) solution to the system (12). The solution is nonnegative for all $i = 1,2,3,...,N$; thus, the summed process (14) also is. The system (12) is often called a continuous-state branching process with multiple types [e.g., 71,72], and is a class of jump-diffusion affine processes [73].

We can directly apply Theorem 2.9 in Barczy et al. [71] to our system (12) for time $t > 0$ if it is conditioned on an initial condition $x_0^{(i)} > 0$ ($i = 1,2,3,...,N$). This follows because the coefficients in the system (12) satisfy the admissibility of parameters required in Definition 2.2 in Barczy et al. [71]. The process $X^{(N)}$ in (14) is therefore well-defined as well. An application of Kolmogorov's two series theorem shows that the limit $N \to +\infty$ of (14) exists if the average and variance of the process $X^{(N)}$ exists as finite values (e.g., Theorem 4.4.4 in Iglói [26]).

Owing to **Assumption 2**, the finite-dimensional model admits a nonnegative average if it exists. Fix $N \in \mathbb{N}$. Assume that the average $m^{(j)}$ of $x^{(j)}$ exists. Then, by taking the expectation of (12), it satisfies

$$r_i m^{(i)} = r_i \frac{b_i}{N} + \frac{M_1}{N}\left(r_i a + \sum_{j=1}^{N} r_j A_{i,j} m^{(j)}\right), \quad i = 1,2,3,...,N , \tag{15}$$

which is rewritten in the vector form

$$\hat{\mathbf{m}} = F(\hat{\mathbf{m}}) \equiv \frac{1}{N}\hat{\mathbf{b}} + \frac{aM_1}{N}\hat{\mathbf{r}} + \frac{1}{N}M_1\mathbf{A}\hat{\mathbf{m}}, \tag{16}$$

where $\hat{\mathbf{m}} = \left[r_i m^{(i)}\right]_{i=1,2,3,\ldots,N}$, $\hat{\mathbf{b}} = \left[r_i b_i\right]_{i=1,2,3,\ldots,N}$, $\hat{\mathbf{r}} = \left[r_i\right]_{i=1,2,3,\ldots,N}$, and $\mathbf{A} = \left[A_{i,j}\right]_{i,j=1,2,3,\ldots,N}$. Due to **Assumption 2**, (16) admits a unique solution $\hat{\mathbf{m}} = \hat{\mathbf{m}}^* \in \mathbb{R}^N$. Indeed, we have the strict contraction

$$\left|F(\hat{\mathbf{m}}_1) - F(\hat{\mathbf{m}}_2)\right| = \frac{1}{N}M_1\left|\mathbf{A}(\hat{\mathbf{m}}_1 - \hat{\mathbf{m}}_2)\right| \le \frac{1}{N}M_1 N\bar{A}\left|\hat{\mathbf{m}}_1 - \hat{\mathbf{m}}_2\right| = \bar{A}M_1\left|\hat{\mathbf{m}}_1 - \hat{\mathbf{m}}_2\right| \tag{17}$$

for any $\hat{\mathbf{m}}_1, \hat{\mathbf{m}}_2 \in \mathbb{R}^N$, where "$|\cdot|$" in (17) is the Euclidean norm in $\mathbb{R}^N$.

We should also check nonnegativity of $\hat{\mathbf{m}}^*$. Each element of $\hat{\mathbf{m}}^*$ is nonnegative if the matrix $\mathbf{B} = I_N - N^{-1}M_1\mathbf{A}$ ($I_N$ is an identity matrix of the dimension $N$) has an inverse whose elements are all nonnegative. It follows that this $\mathbf{B}$ is a non-singular M-matrix whose inverse has nonnegative elements only due to Theorem 2 in Plemmons [74]. Indeed, it is a strictly diagonally dominant Z-matrix because

$$\begin{aligned} B_{i,i} - \sum_{\substack{j=1 \\ i \ne j}}^{N}\left|B_{i,j}\right| &= 1 - \frac{1}{N}M_1 A_{i,i} - \sum_{\substack{j=1 \\ i \ne j}}^{N}\left|B_{i,j}\right| \\ &\ge 1 - \frac{1}{N}\bar{A}M_1 - (N-1)\frac{1}{N}\bar{A}M_1 \quad \text{for } i = 1,2,3,\ldots,N, \\ &= 1 - \bar{A}M_1 \\ &> 0 \end{aligned} \tag{18}$$

which shows $\mathrm{N}_{40}$ of [74] ($\mathbf{B}$ has all positive diagonal elements and it is strictly diagonally dominant) that is equivalent to $\mathrm{F}_{15}$ of the same literature ($\mathbf{B}^{-1}$ exists and only has nonnegative elements. Here, "$D$" in the literature can be selected as an identity matrix). Then, Theorem 1 in Plemmons [74] concludes that each $m^{(i)}$ is nonnegative.

According to Theorem 2.7 of Jin et al. [75], the finite-dimensional system admits a unique stationary probability density if $\nu$ is of the tempered stable type, as assumed later in this study, and the matrix $\mathbf{B}$ is positive definite. The positive-definiteness of $\mathbf{B}$ is satisfied if $A$ is symmetric ($A_{i,j} = A_{j,i}$) because all eigenvalues of a symmetric matrix are real and the condition $\mathrm{C}_9$ in Theorem 1 in Plemmons [74]. A calculation analogous to that in the proof of the Theorem 4.4.4 in Iglói [26] applies at least if $A$ is a sufficiently small constant, and then both the average and variance of $X^{(N)}$ converge to finite values (see the test case in **Section 4**). We consider (and did not prove) that the convergence is in a weak sense (in sense of distribution) because the probability space to define each $\left\{W^{(i)}\right\}_{i=1,2,3,\ldots N}$ and $\left\{J^{(i)}\right\}_{i=1,2,3,\ldots N}$ is different for each $N$ due to the appearance of a wider range of timescales as $N$ grows.

## 2.5 Generalized Riccati equation and statistics

### 2.5.1 Derivation

A generalized Riccati equation is associated with the moment-generating function of an affine process [73], and the moments (hence, cumulants) of the process can be obtained by solving this equation. Moreover, the moment-generating function, which is often used in risk analysis, is an exponential risk indicator for a process [76-79]. We begin with the finite-dimensional version and continue with the original version. The conditional moment-generating function is defined by the exponential moment as follows:

$$\phi_q\left(t,x_1,x_2,...,x_N\right)=\mathbb{E}\left[e^{-qX_t^{(N)}}\middle|\left(x_0^{(1)},x_0^{(2)},...,x_0^{(N)}\right)=\left(x_1,x_2,...,x_N\right)\right],\ \left(x_1,x_2,...,x_N\right)\in\mathbb{R}^N \tag{19}$$

with $q>0$. For each $k\in\mathbb{N}$, we can compute the $k$ th moment as follows:

$$\mathbb{E}\left[\left(X_t^{(N)}\right)^k\middle|\left(x_0^{(1)},x_0^{(2)},...,x_0^{(N)}\right)=\left(x_1,x_2,...,x_N\right)\right]=\left(-1\right)^k\left.\frac{\mathrm{d}^k\phi_q\left(t,x_1,x_2,...,x_N\right)}{\mathrm{d}q^k}\right|_{q\to 0}. \tag{20}$$

Therefore, determining the moment-generating function is crucial for computing the statistics of the nominal model, which can be achieved by solving a generalized Riccati equation.

By exploiting the affine nature of (12), we can infer $\phi_q$ in the exponential form:

$$\phi_q\left(t,x_1,x_2,...,x_N\right)=\exp\left(-\left(F_t+\sum_{i=1}^{N}G_t^{(i)}x_i\right)\right) \tag{21}$$

with time-dependent coefficients $F,G$. According to the definition of $\phi_q$, they satisfy the following initial conditions:

$$F_0=0 \text{ and } G_0^{(i)}=q\ \ (i=1,2,3,...,N). \tag{22}$$

Temporal evolution of $F,G$ is determined through the Kolmogorov equation (e.g., Øksendal and Sulem [14]):

$$\frac{\partial\phi}{\partial t}=\mathfrak{L}\phi,\ \ t>0,\ \left(x_1,x_2,...,x_N\right)\in\mathbb{R}^N \tag{23}$$

with

$$\begin{aligned}\mathfrak{L}\phi&=\sum_{i=1}^{N}\left\{-r_i\left(x_i-\frac{1}{N}b_i\right)\frac{\partial\phi}{\partial x_i}+\frac{1}{2}\sigma^2 r_i x_i\frac{\partial^2\phi}{\partial x_i^2}\right\}+\sum_{i=1}^{N}\frac{1}{N}\left(r_i a+\sum_{j=1}^{N}r_j A_{i,j}x_j\right)\int_0^{+\infty}\Delta\phi\left(z_i\right)v\left(\mathrm{d}z_i\right)\\&=\sum_{i=1}^{N}\left\{\underbrace{-r_i\left(x_i-\frac{1}{N}b_i\right)\frac{\partial\phi}{\partial x_i}}_{\text{Drift}}+\underbrace{\frac{1}{2}\sigma^2 r_i x_i\frac{\partial^2\phi}{\partial x_i^2}}_{\text{Diffusion}}\right\}\\&\quad+\sum_{i=1}^{N}\left\{\underbrace{ar_i\sum_{j=1}^{N}\int_0^{+\infty}\Delta\phi\left(z_j\right)\left(\frac{1}{N}v\left(\mathrm{d}z_j\right)\right)}_{\text{Exogenous jump}}+\underbrace{x_i r_i\int_0^{+\infty}\Delta\phi\left(z_j\right)\left(\sum_{j=1}^{N}\frac{1}{N}A_{j,i}v\left(\mathrm{d}z_j\right)\right)}_{\text{Self-exciting jump}}\right\}\end{aligned} \tag{24}$$

and

$$\Delta\phi\left(z_i\right)=\left(x_1,x_2,...,x_i+z_i,...,x_N\right)-\phi. \tag{25}$$

Substituting (21) into (23) yields

$$\frac{\mathrm{d}}{\mathrm{d}t}\left(-F_t-\sum_{i=1}^{N}G_t^{(i)}x_i\right)=\sum_{i=1}^{N}\left\{-r_i\left(x_i-\frac{1}{N}b_i\right)\left(-G_t^{(i)}\right)+\frac{1}{2}\sigma^2 r_i x_i\left(-G_t^{(i)}\right)^2\right\} \\ +\sum_{i=1}^{N}\left\{\frac{1}{N}r_i a\int_0^{+\infty}\left(e^{-G_t^{(i)}z_i}-1\right)v\left(\mathrm{d}z_i\right)+\frac{r_i}{N}x_i\sum_{j=1}^{N}A_{j,i}\int_0^{+\infty}\left(e^{-G_t^{(j)}z_i}-1\right)v\left(\mathrm{d}z_j\right)\right\}, \tag{26}$$

from which we obtain the generalized Riccati equation as follows:

$$\frac{\mathrm{d}F_t}{\mathrm{d}t}=\frac{1}{N}\sum_{i=1}^{N}r_i\left\{b_iG_t^{(i)}+a\int_0^{+\infty}\left(1-e^{-G_t^{(i)}z_i}\right)v\left(\mathrm{d}z_i\right)\right\} \tag{27}$$

and

$$\frac{\mathrm{d}G_t^{(i)}}{\mathrm{d}t}=-r_iG_t^{(i)}-\frac{1}{2}\sigma^2 r_i\left(G_t^{(i)}\right)^2+r_i\frac{1}{N}\sum_{j=1}^{N}A_{j,i}\int_0^{+\infty}\left(1-e^{-G_t^{(j)}z_j}\right)v\left(\mathrm{d}z_j\right) \tag{28}$$

subject to the initial conditions (22).

The generalized Riccati equation associated with the nominal model is derived by formally letting $N\to+\infty$ in (27)–(28) and (22) with some abuse of notations:

$$\frac{\mathrm{d}F(t)}{\mathrm{d}t}=\int_{r=0}^{r=+\infty}r\left\{b(r)G(t,r)+a\int_{z=0}^{z=+\infty}\left(1-e^{-G(t,r)z}\right)v(\mathrm{d}z)\right\}\rho(\mathrm{d}r),\quad t>0 \tag{29}$$

and

$$\frac{\mathrm{d}G(t,r)}{\mathrm{d}t}=-rG(t,r)-\frac{1}{2}\sigma^2 r\left(G(t,r)\right)^2 \\ +r\int_{u=0}^{u=+\infty}A(u,r)\int_{z=0}^{z=+\infty}\left(1-e^{-G(t,u)z}\right)v(\mathrm{d}z)\rho(\mathrm{d}u),\quad t>0,\ r>0 \tag{30}$$

subject to initial conditions

$$F(0)=0\ \text{ and }\ G(0,\cdot)=q, \tag{31}$$

which is a system of nonlinear partial integro-differential equations.

***Remark 1*** One may also consider the exponential moment with expectant $e^{qX_t^{(N)}}$ with $q>0$ to evaluate the sensitivity of the upper tail of $X$. In this case, we should *a priori* assume that $q$ is sufficiently close to zero; otherwise, the corresponding exponential moment may diverge in finite time, as suggested for a class of affine jump-diffusion processes [80]. The critical value of $q$ for ensuring the global existence of a moment has been determined in some cases [81,82]. See **Section 2.5.3**.

### 2.5.2 Mathematical analysis of generalized Riccati equation

Generalized Riccati equations in both finite- and infinite-dimensional forms are studied by focusing on the unique existence of solutions. We focus on the infinite-dimensional case because the results for the finite-dimensional case followed the formal replacement $\rho\to\rho_N$.

The following technical lemma is used to investigate the generalized Riccati equation.

***Lemma 1***

*For* $q \in \left(0, \frac{2}{\sigma^2}\left(1-\bar{A}M_1\right)\right)$*, there exists a unique solution* $\varphi = \bar{\varphi}(q) > 0$ *to the following equation:*

$$\varphi = q + \frac{\bar{A}M_1 + \frac{\sigma^2}{2}\varphi}{1+\frac{\sigma^2}{2}\varphi}\varphi \ \left(:= \gamma(\varphi)\right), \quad \varphi \geq 0. \tag{32}$$

*Moreover,* $\bar{\varphi}(q) > q$ *and* $\lim_{q \to +0} \bar{\varphi}(q) = 0$.

We set the space of the bounded continuous functions equipped with the supremum norm $\|\cdot\|$, which is a Banach space (e.g., p.30 in Clason [83]):

$$C_{\mathrm{b}} = \left\{ \varphi \in C(0,+\infty) : \|\varphi\| = \sup_{r \in (0,+\infty)} |\varphi(r)| < +\infty \right\}. \tag{33}$$

For each $T > 0$, we set another Banach space $C_{\mathrm{b},T}$ of bounded functions as follows:

$$C_{\mathrm{b},T} = \left\{ \varphi \in C\left([0,T]\times(0,+\infty)\right) : \|\varphi\|_T = \sup_{0 \leq t \leq T} \|\varphi(t,\cdot)\| < +\infty \right\}. \tag{34}$$

We set $\eta(q) = 1 + \frac{\sigma^2}{2}\bar{\varphi}(q)$.

By a solution to the second generalized Riccati equation (30), we mean $\varphi \in C_{\mathrm{b},T}$ obtained from the variation of constants formula

$$\begin{aligned} \varphi(t,r) &= q e^{-\eta(q)rt} \\ &\quad + \int_{s=0}^{s=t} e^{-\eta(q)r(t-s)} \frac{1}{2} r\sigma^2 \left(\bar{\varphi}(q) - \varphi(s,r)\right)\varphi(s,r)\,\mathrm{d}s \\ &\quad + \int_{s=0}^{s=t} e^{-\eta(q)r(t-s)} r \int_{u=0}^{u=+\infty} A(u,r) \int_{z=0}^{z=+\infty} \left(1 - e^{-\varphi(s,u)z}\right) v(\mathrm{d}z)\rho(\mathrm{d}u)\,\mathrm{d}s \\ &\quad \left(:= \mathbb{G}(\varphi)(t,r)\right) \end{aligned} \tag{35}$$

for $0 \leq t \leq T$ and $r > 0$, for each $T > 0$. The functional form of (35) is derived from the formal representation of (30):

$$\begin{aligned} \frac{\mathrm{d}}{\mathrm{d}t}\left(e^{\eta(q)rt} G(t,r)\right) &= e^{\eta(q)rt} \frac{1}{2} r\sigma^2 \left(\bar{\varphi}(q) - G(t,r)\right) G(t,r) \\ &\quad + e^{\eta(q)rt} r \int_{u=0}^{u=+\infty} A(u,r) \int_{z=0}^{z=+\infty} \left(1 - e^{-G(t,u)z}\right) v(\mathrm{d}z)\rho(\mathrm{d}u) \end{aligned}. \tag{36}$$

We obtain a proposition regarding the unique existence of a solution to (30) in the sense of exists uniquely.

***Proposition 1***

*For* $q \in \left(0, \frac{2}{\sigma^2}\left(1-\bar{A}M_1\right)\right)$, *the second generalized Riccati equation (30) admits a unique nonnegative solution that is bounded as* $0 \le \varphi(t,r) \le \bar{\varphi}(q)$ *for all* $t \ge 0$ *and* $r > 0$. *Moreover, the solution is continuously differentiable at any time* $t > 0$, *and thus satisfies (30) pointwise.*

Considering **Proposition 1**, we obtain the following result concerning the first Riccati equation (29). A solution to the first Riccati equation (29) is a bounded function $\psi \in C([0,+\infty))$ such that

$$\psi(t) = \int_0^t \int_{r=0}^{r=+\infty} r\left\{b(r)G(s,r) + a\int_{z=0}^{z=+\infty}\left(1-e^{-G(s,r)z}\right)v(\mathrm{d}z)\right\}\rho(\mathrm{d}r)\mathrm{d}s, \quad t \ge 0 \tag{37}$$

with $G$ being the unique nonnegative and bounded solution to (30).

***Proposition 2***

*For a small* $q > 0$, *the first Riccati equation admits a unique solution that is nonnegative.*

### 2.5.3 Positive exponential case

It is also possible to consider a generalized Riccati equation corresponding to the moment of $e^{qX_t}$ ($q > 0$). In this case, the equations to determine the moment, with the abuse of notation, are as follows:

$$\frac{\mathrm{d}F(t)}{\mathrm{d}t} = \int_{r=0}^{r=+\infty} r\left\{b(r)G(t,r) + a\int_{z=0}^{z=+\infty}\left(e^{G(t,r)z}-1\right)v(\mathrm{d}z)\right\}\rho(\mathrm{d}r), \quad t > 0 \tag{38}$$

and

$$\begin{aligned}\frac{\mathrm{d}G(t,r)}{\mathrm{d}t} &= -rG(t,r) + \frac{1}{2}\sigma^2 r\left(G(t,r)\right)^2 \\ &\quad + r\int_{u=0}^{u=+\infty} A(u,r)\int_{z=0}^{z=+\infty}\left(e^{G(t,u)z}-1\right)v(\mathrm{d}z)\rho(\mathrm{d}u)\end{aligned}, \quad t > 0, \quad r > 0 \tag{39}$$

subject to the initial conditions (31). A solution to (39) is defined through

$$\varphi(t,r) = qe^{-rt} + \int_{s=0}^{s=t} e^{-r(t-s)} r\left\{\frac{1}{2}\sigma^2\left(\varphi(s,r)\right)^2 + \int_{u=0}^{u=+\infty} A(u,r)\int_{z=0}^{z=+\infty}\left(e^{\varphi(s,u)z}-1\right)v(\mathrm{d}z)\rho(\mathrm{d}u)\right\}\mathrm{d}s. \tag{40}$$

An elementary calculation analogous to **Lemma 1** shows that the equation

$$\varphi = q + \frac{1}{2}\sigma^2\varphi^2 + \bar{A}\int_{z=0}^{z=+\infty}\left(e^{\varphi z}-1\right)v(\mathrm{d}z) \ \left(:= \vartheta(\varphi)\right), \quad \varphi > 0 \tag{41}$$

admits a positive solution, again denoted by $\bar{\varphi}(q)$, if

$$\frac{\mathrm{d}\vartheta(\varphi)}{\mathrm{d}\varphi} = \sigma^2\bar{\varphi}(q) + \bar{A}\int_{z=0}^{z=+\infty} ze^{\bar{\varphi}(q)z}v(\mathrm{d}z) < 1 \ \text{ and } \ \vartheta(\varphi) > \varphi \ \text{ in } \ (0,\xi) \tag{42}$$

for a small constant $\xi > 0$, such that $\int_{z=0}^{z=+\infty}\left(e^{\xi z}-1\right)v(\mathrm{d}z) < +\infty$ and $\int_{z=0}^{z=+\infty} ze^{\xi z}v(\mathrm{d}z) < 1$, for $q \in (0,\xi)$ chosen to be sufficiently small. We consider the following auxiliary equation:

$$\varphi(t,r) = qe^{-rt} + \int_{s=0}^{s=t} e^{-r(t-s)} r \left\{ \frac{1}{2}\sigma^2 \left(\hat{\varphi}(s,r)\right)^2 + \int_{u=0}^{u=+\infty} A(u,r) \int_{z=0}^{z=+\infty} \left(e^{\hat{\varphi}(s,u)z} - 1\right) \nu(\mathrm{d}z) \rho(\mathrm{d}u) \right\} \mathrm{d}s \\ \left(:= \hat{\mathbb{K}}(\varphi)(t,r)\right) \tag{43}$$

with the truncation $\hat{\varphi}$ using the new $\bar{\varphi}(q)$. Based on the strategy analogous to **Proof of Proposition 1**, there exists a unique bounded continuous solution $0 \le \tilde{\varphi} \le \bar{\varphi}(q)$, which is also a solution to (39) if $q$ is sufficiently small. This is owing to the continuity of the right-hand side of (43) along with the boundedness

$$\left\| \hat{\mathbb{K}}(\varphi_1) \right\|_T \le \bar{\varphi}(q) \tag{44}$$

and the strict contraction (we use (42))

$$\left\| \hat{\mathbb{K}}(\varphi_1) - \hat{\mathbb{K}}(\varphi_2) \right\|_T \le \left( \sigma^2 \bar{\varphi}(q) + \bar{A} \int_{z=0}^{z=+\infty} z e^{\bar{\varphi}(q)z} \nu(\mathrm{d}z) \right) \left\| \varphi_1 - \varphi_2 \right\|_T \tag{45}$$

for any $\varphi_1, \varphi_2 \in C_{\mathrm{b},T}$, for each $T > 0$. The uniqueness of the solutions to (39) follows from (42) and

$$\left\| \varphi_1 - \varphi_2 \right\|_T \le \left( \sigma^2 \bar{\varphi}(q) + \bar{A} \int_{z=0}^{z=+\infty} z e^{\bar{\varphi}(q)z} \nu(\mathrm{d}z) \right) \left\| \varphi_1 - \varphi_2 \right\|_T . \tag{46}$$

## 3. Numerical scheme

### 3.1 Exact discretization methods

The exact discretization methods considered in this study are the one for square-root diffusions (48) and that for tempered stable Ornstein–Uhlenbeck processes (54). Their detailed explanations are available in the corresponding references cited in the sequel. For modeling simplicity, we assume the tempered-stable type $\nu$ serving as a principal model of driving jump processes [84-87]:

$$\nu(\mathrm{d}z) = \frac{a_\nu e^{-b_\nu z}}{z^{1+c_\nu}} \mathrm{d}z , \quad z > 0 \tag{47}$$

with $a_\nu > 0$ that controls jump frequency, $b_\nu > 0$ that tilts jump sizes, and $c_\nu < 1$ that controls fluctuation of small jumps. The tempered-stable model (47) serves as a versatile mathematical tool that generates jumps with positive and bounded variations with ( $c_\nu < 0$ , compound Poisson) and without ( $0 \le c_\nu < 1$, not compound Poisson) boundedness of jump events in each fixed time interval.

#### 3.1.1 Exact discretization of square-root diffusion

The method proposed by Abi Jaber [57] discretizes the following square-root diffusion for any $a_1, b_1, c_1 \in \mathbb{R}$:

$$\mathrm{d}Z_t = \left(a_1 + b_1 Z_t\right)\mathrm{d}t + c_1 \sqrt{Z_t}\,\mathrm{d}W_t , \quad t > 0 \tag{48}$$

subject to $Z_0 \ge 0$, where $W$ is a one-dimensional standard Brownian motion. As implied in the nominal model, the cases of interest correspond to $a_1 \ge 0$ and $b_1 \le 0$.

The SDE (48) admits a unique strong solution whose boundary behavior critically depends on the parameter values, and the solution is always positive for $Z_0 > 0$ if the diffusion is small ($a_1 \geq 2c_1^2$), indicating that the so-called Feller condition is satisfied, while the solution touches boundary 0 in a finite time with probability one if the diffusion is large ($a_1 < 2c_1^2$) (Proposition 1.2.15 in Alfonsi [59]). Intriguingly, the latter with a large diffusion, is a more difficult subject because the non-Lipschitz nature of the diffusion coefficient $\sqrt{Z_t}$ becomes dominant.

The method proposed by Abi Jaber [57] directly discretizes the time integration $\int_0^t Z_s \mathrm{d}s$, and subsequently obtains a discretization scheme of the process $Z$ in a way that the discretized version of (48) corresponds to the associated Riccati equation. The resulting scheme unconditionally preserves nonnegativity of numerical solutions to the SDE (48) and provides accurate results against several test cases as well as economic [57] and environmental case studies [40]. The scheme is expressed as follows, where $\Delta t$ represents time increment and $Z_n$ represents the approximation of $Z$ at time $n\Delta t$ ($n = 0,1,2,\ldots$):

$$Z_{n+1} = Z_n + a_1 \Delta t + b_1 U_n + c_1 V_n, \quad n = 0,1,2,\ldots, \tag{49}$$

where, $U_n$ and $V_n$ in (49) are expressed as follows: if $b_1 > 0$,

$$U_n \sim \mathrm{InvGamma}\left(\kappa_n, \left(\frac{\kappa_n}{\omega_n}\right)^2\right) \text{ and } V_n = \frac{1}{\omega_n}\left(U_n - \kappa_n\right) \tag{50}$$

with

$$\kappa_n = Z_n \frac{e^{b_1 \Delta t} - 1}{b_1} + \frac{a_1}{b_1}\left(\frac{e^{b_1 \Delta t} - 1}{b_1} - \Delta t\right) \text{ and } \omega_n = c_1 \frac{e^{b_1 \Delta t} - 1}{b_1}. \tag{51}$$

If $b_1 = 0$, $\kappa_n$ and $\omega_n$ become

$$\kappa_n = Z_n \Delta t + \frac{1}{2} a_1 \left(\Delta t\right)^2 \text{ and } \omega_n = c_1 \Delta t. \tag{52}$$

Here, the inverse Gamma distribution $\mathrm{InvGamma}\left(\mu, \lambda\right)$ admits the following density:

$$p_{IG}\left(y\right) = \sqrt{\frac{\lambda}{2\pi y^2}} \exp\left(-\frac{\lambda\left(y-\mu\right)^2}{2\mu^2 y}\right), \quad y > 0. \tag{53}$$

Each $U_n$ is assumed to be mutually independent. Therefore, the scheme discretizes the SDE (48) as the sum of deterministic and random increments, where the inverse gamma noise replaces the Gaussian noise in the common Euler–Maruyama method. The name "exact discretization" comes from the fact that in theory the equation (49) correctly reproduces the average of the process $Z$ for any $\Delta t$. We do not use the Euler–Maruyama method because it works poorly if $c_1$ is large [e.g. 56], while the method proposed by Abi Jaber [57] works for the whole $c_1 > 0$.

***Remark 2*** The case $b_1 = 0$ corresponds to a Bessel process to which a tailored numerical scheme applies [88]. We use the method examined above because it applies to both $b_1 > 0$ and $b_1 = 0$.

### 3.1.2 Exact discretization of tempered stable Ornstein–Uhlenbeck process

Another exact discretization method used in this study involves discretizing the following jump-driven Ornstein–Uhlenbeck process for any $a_2 > 0$:

$$\mathrm{d}S_t = -a_2 S_t \mathrm{d}t + \mathrm{d}L_t, \quad t > 0 \tag{54}$$

subject to $S_0 \geq 0$, where $L$ is a one-dimensional pure-jump Lévy process with a Lévy measure (47).

The numerical method proposed by Qu et al. [58] discretizes the SDE (54) based on the consistency between the conditional moment-generating function and the discretized SDE; thus, numerical solutions are unconditionally nonnegative, and both large and small jumps generated by the driving Lévy process $L$ are accounted for at each time step. The full explanation of the exact discretization method is not presented in this paper because it requires some nested computation at each time step with composite coefficients and random variables; however, these random variables are computable using acceptance rejection methods (Algorithms 3.1 and 3.2 in Qu et al. [58]).

Using the notations in the previous subsection, the scheme is expressed as follows:

$$S_{n+1} = S_n e^{-a_2 \Delta t} + \Theta_n + \sum_{k=1}^{\Lambda_n} \Xi_{k,n}, \quad n = 0, 1, 2, \ldots, \tag{55}$$

where $\Theta_n$ is generated from a specific tempered stable distribution (Algorithm D.1 in Qu et al. [58]), $\Xi_{k,n}$ is generated from a specific gamma distribution whose scaling parameter is based on another random variable (Algorithm 3.2 in Qu et al. [58]), and $\Lambda_n$ is a Poisson random variable with intensity:

$$\frac{a_v}{c_v a_2} b_v^{c_v} \Gamma(1 - c_v) \left( \frac{1}{c_v} \left( \left( e^{-a_2 \Delta t} \right)^{-c_v} - 1 \right) - \ln\left( e^{-a_2 \Delta t} \right) \right). \tag{56}$$

The name "exact discretization" comes from the fact that the right-hand side of (55) equals the expectation of $S_{n+1}$ conditioned on $S_n$.

***Remark 3*** Some studies were conducted on Ornstein–Uhlenbeck processes driven by other jump processes [89-91]. They may replace the method proposed by Qu et al. [58] when necessary.

***Remark 4*** Arai and Imai [92] highlights some uncertainties in the method proposed by Qu et al. [58]. This does not affect our case because it corresponds to "$\rho = 1$" in their notation.

## 3.2 Proposed scheme

The proposed scheme for the nominal model combines two discretization methods. We consider the following three cases: (a) no jump, (b) no diffusion, and (c) both jumps and diffusions. Cases (a) and (b) use the methods proposed by Abi Jaber [57] and Qu et al. [58], respectively. Case (c) uses both methods, along with operator splitting.

### 3.2.1 Case (a): no jump

In this case, the nominal model (12) discretized for $r > 0$ reads

$$\mathrm{d}x_t^{(i)} = -r_i x_t^{(i)}\mathrm{d}t + \frac{1}{N} r_i b_i \mathrm{d}t + \sigma\sqrt{r_i x_t^{(i)}}\mathrm{d}W_t^{(i)}, \quad t \in \mathbb{R}, \quad i = 1,2,3,...,N. \tag{57}$$

This SDE is a reparameterization of (48) with the replacements $a_1 \to \frac{1}{N} r_i b_i$, $b_1 \to -r_i$, and $c_1 \to \sigma\sqrt{r_i}$. Therefore, we apply the method proposed by Abi Jaber [57] to each $i$.

### 3.2.2 Case (b): no diffusion

In this case, the nominal model (12) discretized for $r > 0$ reads

$$\mathrm{d}x_t^{(i)} = -r_i x_t^{(i)}\mathrm{d}t + \int_{u=0}^{u=r\frac{1}{N}(a+\mu_i)} \int_{z=0}^{z=+\infty} zJ^{(i)}(\mathrm{d}t, \mathrm{d}u, \mathrm{d}z), \quad t \in \mathbb{R}, \quad i = 1,2,3,...,N. \tag{58}$$

This SDE is just a reparameterization of the Ornstein–Uhlenbeck process (54). Subsequently, we applied the method proposed by Qu et al. [58] with a fully explicit treatment of the state-dependent coefficient $\mu_i$. We discretize $\mu_i$ at time $n\Delta t$ as follows:

$$\mu_i \approx \mu_{i,n} = \sum_{j=1}^{N} A_{i,j} x_{n\Delta t}^{(j)}. \tag{59}$$

The state $x^{(j)}$ in the jump term is thus frozen at time $n\Delta t$; thus, $x_{n\Delta t}^{(j)}$ is treated as a known quantity while computing $x_{(n+1)\Delta t}^{(j)}$. Subsequently, the numerical solution at time $t = (n+1)\Delta t$ is obtained directly using the method proposed by Qu et al. [58], along with the following replacement: $a_2 \to r_i$ and $a_v \to a_v r \frac{1}{N}(a + \mu_{i,n})$.

### 3.2.3 Case (c): both jump and diffusion

We apply the operator-splitting method to (12). Assuming that all $x_{n\Delta t}^{(i)}$ ($i = 1,2,3,...,N$) are available, we use the following **Algorithm 1**.

***Algorithm 1***

***Step 1.*** *Discretize the SDE using the method proposed by Abi Jaber [57] in the time interval* $(n\Delta t, (n+1)\Delta t)$. *The updated numerical solution is denoted by* $\hat{x}_{n\Delta t}^{(i)}$ ($i = 1,2,3,...,N$)*:*

$$\mathrm{d}x_t^{(i)} = \frac{1}{N} r_i b_i \mathrm{d}t + \sigma\sqrt{r_i x_t^{(i)}}\mathrm{d}W_t^{(i)}, \quad i = 1,2,3,...,N. \tag{60}$$

***Step 2.*** *Discretize the SDE (58) using the method proposed by Qu et al. [58] with the time increment* $\Delta t$ *in the time interval* $(n\Delta t,(n+1)\Delta t)$ *starting from the initial value* $\hat{x}_{n\Delta t}^{(i)}$ *(*$i=1,2,3,...,N$*). The updated numerical solution is denoted by* $x_{(n+1)\Delta t}^{(i)}$ *(*$i=1,2,3,...,N$*).*

A theoretical advantage of **Algorithm 1** is that it generates nonnegative superposition processes because each $x_t^{(i)}$ at each step is necessarily nonnegative, and the algorithm is an alternating iteration for the diffusion and jump parts. Applying the Euler–Maruyama method never associates such a strong property, particularly for Step 1, because it has a nonzero probability of obtaining a negative $x_t^{(i)}$ value owing to the unboundedness of the discretized Brownian motion $W_t^{(i)}$.

## 4. Computational investigations

This section has the following two goals. The first objective is to evaluate the accuracy of the proposed discretization method for the nominal model. The second goal is to examine the numerical method against more complex cases: a water quantity-quality model and a model with nonlinearity that is not affine. The role of **Section 4.1** is to validate the proposed numerical scheme against nonlinearity, particularly the splitting method because it deals with the extreme square-root diffusion case where the drift completely vanishes, which cannot be dealt with by classical numerical schemes such as the Euler–Maruyama method [e.g., 56].

### 4.1 Numerical tests against the nominal model

The first test case considers the square-root SDE

$$\mathrm{d}Z_t = c_1\sqrt{Z_t}\,\mathrm{d}W_t\,,\quad t>0 \tag{61}$$

under $Z_0>0$. This SDE is a building block of the operator-splitting method. Without loss of generality, we set $c_1=1$. We compute the sample paths using the exact discretization method and evaluated their average and variance, which are non-pathwise quantities, to test the weak approximation ability (i.e., moment approximation ability) of the scheme and the following first hitting time $\tau$ as a pathwise quantity:

$$\tau = \inf\left\{t>0 \middle| Z_t = 0 \text{ or } 1\right\}. \tag{62}$$

We then compute the first-order moment $T_1$ and second-order moment $T_2$ of $\tau$. The average and variance of $Z_t$ are obtained directly from (61):

$$\mathbb{E}[Z_t]=Z_0 \text{ and } \mathbb{V}[Z_t]=\sigma^2 Z_0 t\,,\quad t\geq 0\,. \tag{63}$$

The moments $T_1$ and $T_2$ are derived by solving the following Kolmogorov equations [e.g., 93]:

$$\frac{1}{2}Z_0\frac{\mathrm{d}^2T_k}{\mathrm{d}Z_0^2} = -kT_{k-1} \text{ for } 0<Z_0<1\,,\ T_k(0)=T_k(1)=0 \tag{64}$$

for $k = 1, 2$, where $T_0 = 1$. These equations are solved for $0 < Z_0 \le 1$ as follows:

$$\mathbb{E}[\tau] = T_1 = -2Z_0 \ln Z_0 \tag{65}$$

and

$$\mathbb{V}[\tau] = T_2 - T_1^2 = 4Z_0^2 \ln Z_0 - 4Z_0^2 (\ln Z_0)^2 + 6Z_0 (1 - Z_0). \tag{66}$$

The total simulation duration for each sample path is 10 that is sufficiently long for our purpose. We examine the numerical method for different values of the time increment $\Delta t$ and the total number of sample paths $N_{\mathrm{P}}$. The error in this test case indicates an absolute error. For the average and variance of $Z$, the absolute error is the maximum value of the difference between the theoretical and computed results at all time steps. For the average and variance of $\tau$, the error is the difference between the theoretical and computed results.

**Figure 1** shows the sample paths of $Z$, and **Figure 2** compares the average and variance of $Z$ with $Z_0 = 0.5$, $\Delta t = 0.0001$, and $N_{\mathrm{P}} = 100,000$ paths. **Table 1** lists a comparison of the absolute differences between the theoretical and computed averages. **Table 2** lists the corresponding results for variance. **Table 3** lists the errors in the computed averages of $\tau$. Similarly, **Table 4** lists the errors in the computed variance of $\tau$. We must balance the time increment and the total number of sample paths, as suggested in **Tables 1–4** because only using a sufficiently small $\Delta t$ results in an increase in the error, for example comparing cases $\Delta t = 0.000050$ and $\Delta t = 0.000025$ with $N_{\mathrm{P}} = 100,000$. For the computational cases presented in this subsection, it appears reasonable to set $(\Delta t, N_{\mathrm{P}}) = (1/2^l \times 0.0001,\ 2^l \times 100,000)$ ($l = 0,1,2$), which are the diagonals in **Tables 1–4**, with which the error decreases as $l$ increases. This conjecture is supported by the additional computational results for $l = 3$; the errors in average and variance of $Z$ are 2.647.E-03 and 3.978.E-02, respectively, and the errors in average and variance of $\tau$ are 3.504.E-03 and 5.635.E-03, respectively. The errors in the average and variance of $\tau$ decreases as expected and were almost half of those with $l = 2$, whereas those of $Z$ do not at this resolution. Further doubling $N_{\mathrm{P}}$ while $\Delta t$ is fixed yields smaller errors in the average and variance of $Z$, which are 1.892.E-03 and 1.947.E-02, respectively. These observations, provided that the error with respect to $N_{\mathrm{P}}$ scales $N_{\mathrm{P}}^{-0.5}$ as usual in Monte Carlo simulations, suggest that the convergence rate of the discretization method in time had an order larger than 0.5, for cumulants of $Z$, at least for the cases examined here. Although a detailed theoretical convergence study of the numerical method should be conducted in the future, the computational results obtained suggest that the numerical method proposed by Abi Jaber [57] could potentially approximate the SDE (61) in both weak and strong senses if the resolution in time and the total number of samples were balanced. This implies that the nonlinearity in the nominal model can be handled by the proposed numerical scheme.

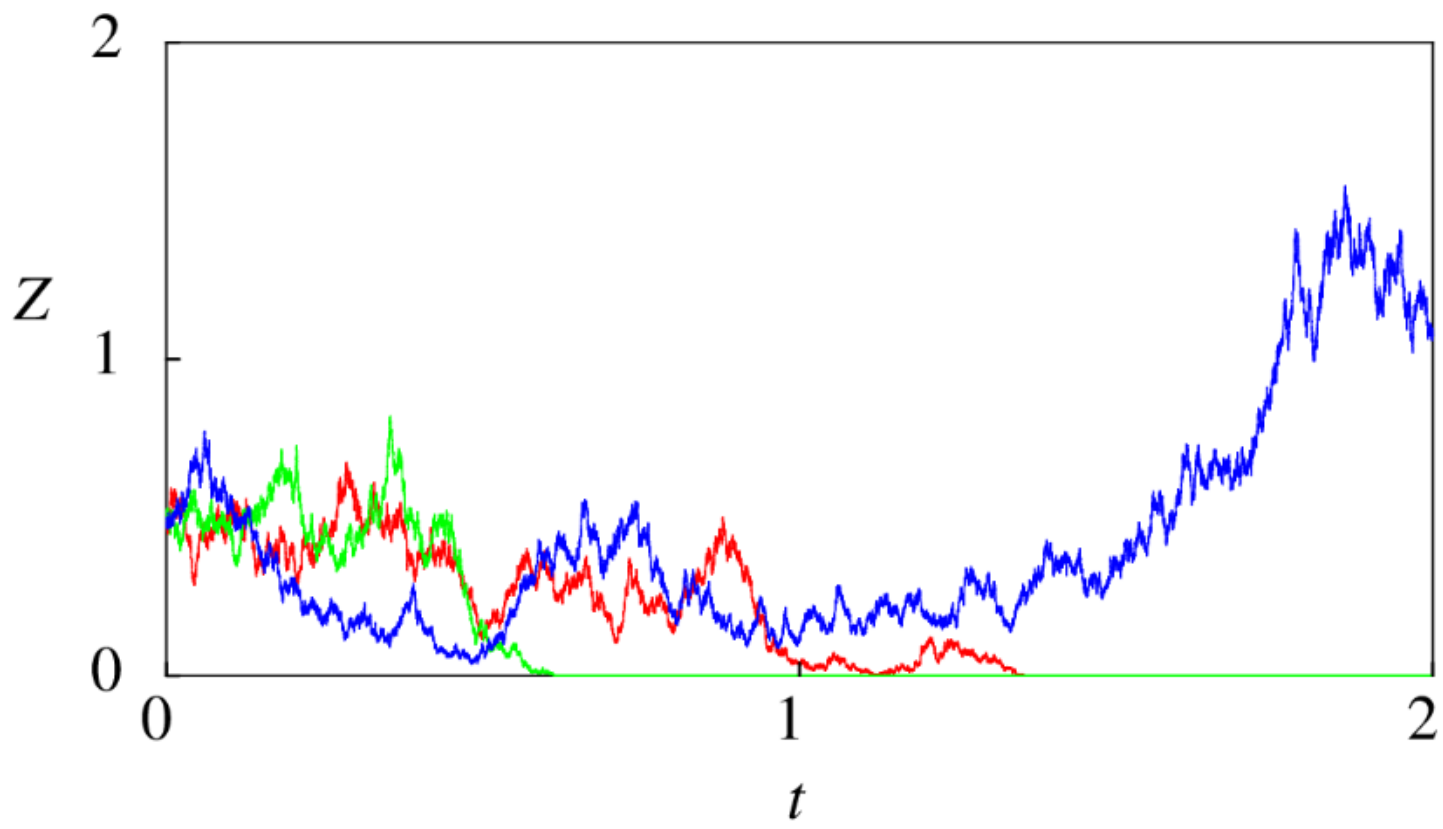

**Figure 1.** Computed sample paths of the square-root SDE. Distinct colors represent different sample paths.

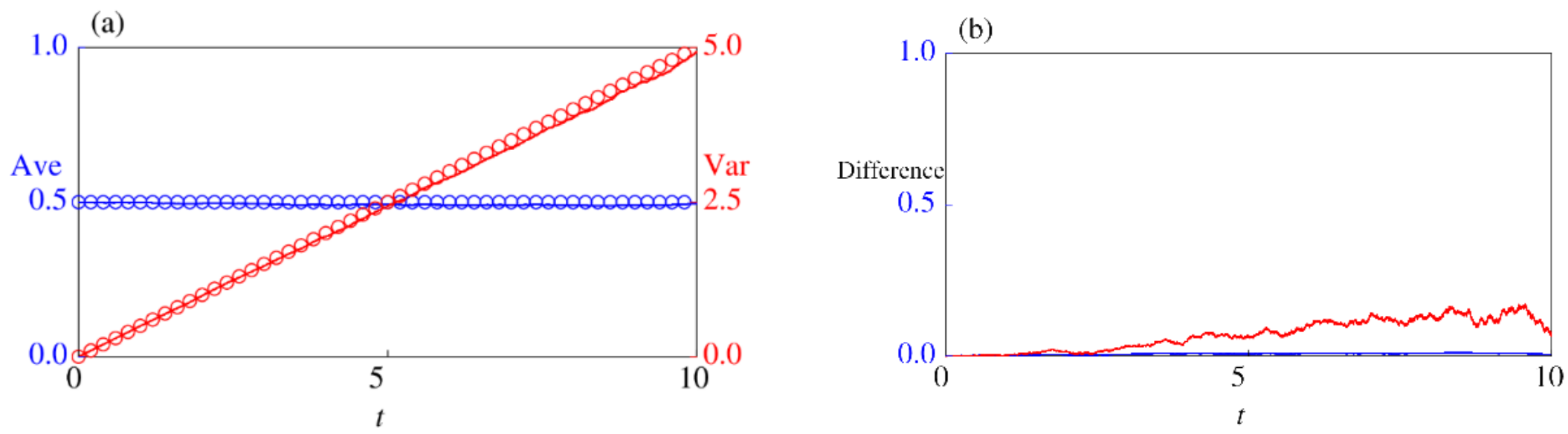

**Figure 2.** (a) Comparison of theoretical and computed average (red) and variance (blue): Circles (theoretical results) and Lines (computational results) and (b) their difference with the same color legends.

**Table 1.** Relative error between theoretical and computed average of $Z$.

| | | $N_P$ | | |
|---|---|---|---|---|
| | | 100,000 | 200,000 | 400,000 |
| | 0.000100 | 1.227.E-02 | 9.369.E-03 | 2.780.E-03 |
| $\Delta t$ | 0.000050 | 8.273.E-03 | 5.552.E-03 | 3.650.E-03 |
| | 0.000025 | 6.268.E-03 | 2.568.E-03 | 2.844.E-03 |

**Table 2.** Relative error between theoretical and computed variance of $Z$.

| | | $N_P$ | | |
|---|---|---|---|---|
| | | 100,000 | 200,000 | 400,000 |
| | 0.000100 | 1.713.E-01 | 1.633.E-01 | 5.545.E-02 |
| $\Delta t$ | 0.000050 | 1.112.E-01 | 7.982.E-02 | 9.734.E-02 |
| | 0.000025 | 1.061.E-01 | 7.868.E-02 | 2.972.E-02 |

**Table 3.** Relative error between theoretical and computed average $\tau$. The theoretical value is 0.69315.

| | | $N_P$ | | |
|---|---|---|---|---|
| | | 100,000 | 200,000 | 400,000 |
| | 0.000100 | 1.617.E-02 | 1.449.E-02 | 1.240.E-02 |
| $\Delta t$ | 0.000050 | 9.249.E-03 | 6.317.E-03 | 7.252.E-03 |
| | 0.000025 | 7.304.E-03 | 5.883.E-03 | 5.239.E-03 |

**Table 4.** Relative error between theoretical and computed variance of $\tau$. The theoretical value is 0.32640.

| | | $N_P$ | | |
|---|---|---|---|---|
| | | 100,000 | 200,000 | 400,000 |
| | 0.000100 | 2.168.E-02 | 2.927.E-02 | 2.382.E-02 |
| $\Delta t$ | 0.000050 | 2.105.E-02 | 1.631.E-02 | 1.643.E-02 |
| | 0.000025 | 2.828.E-02 | 1.673.E-02 | 1.081.E-02 |

### 4.2 Nominal model

#### 4.2.1 Computational conditions

The investigations in the remaining part of **Section 4** focus on the ability of the proposed numerical method to deal with the superposition processes against different values of the degree of freedom $N$. We assume that all parameters and variables are non-dimensional for simplicity. We also assume that the probability measure $\rho$ is a gamma distribution with the shape parameter $\alpha$ and scale parameter $\beta$ (i.e., average $\alpha\beta$ and variance $\alpha\beta^2$) that have been often assumed in the superposition approaches [e.g., 39, 67]. Throughout the numerical test cases, we fix the parameter value $c_v = 0.65$ because in applications the range of $c_v$ approximately 0.5 to 0.7 appeared to be realistic in engineering applications [e.g., 27, 69]. In both Tests 1 and 2, we set $\sigma = 0.5$ and fixed $b_v = 0.1$ and chose $a_v$ so that $M_1 = 1/2$ for simplicity. Statistics and probability density are computed by assuming ergodicity; thus, an ensemble average is replaced by a time average. Selecting the proper length of the burn-in period for each stochastic process is in general a difficult task, and we examine it empirically. We consider that the risk that the ergodic assumption is broken is small considering the convergence result (e.g., **Figure 4** and **Table 5**) and the mean-reverting nature of the target process. We sample statistics for the length of time for 400,000 or longer one 1,600,000. The initial condition is $x_0^{(i)} = 0$ for all $i$. Unless otherwise specified, we fix $\Delta t = 0.01$.

We consider two cases of the self-excitation $\mu_x$ of jumps. Both cases use

$$\mu_x(r) = rA_0 x_t(\mathrm{d}r) + \int_0^{+\infty} uA_1(r,u)x_t(\mathrm{d}u), \tag{67}$$

which is a variant of (9). The first case use $A_0 \equiv 1$ and $A_1 \equiv 0$ (Test 1), while the second use $A_0 \equiv 0$ and $A_1 \equiv 1$ (Test 2); thus, the second case assume self-excitation owing to aggregation. The average and variance of the first case in the stationary state are obtained through formal calculations as follows:

$$\mathbb{E}[X_t] = \frac{M_1}{1-M_1}a \quad \text{and} \quad \mathbb{V}[X_t] = \frac{a}{2(1-M_1)^2}\left(\sigma^2 M_1 + M_2\right). \tag{68}$$

The average is free of $\sigma$, whereas the variance is not. The autocorrelation function with lag $l \geq 0$ is obtained as $(1 + M_1\beta l)^{-\alpha}$, and can be compared with numerical functions. For Test 1, we consider both moderate ($\alpha = 2.5$, and thus $\alpha > 1$) and long memory ($\alpha = 0.8$, and thus $0 < \alpha < 1$). We fix $\beta = 1$.

The cumulants of Test 2 have not been found explicitly; however, the average is obtained by solving the corresponding linear system, as described in **Section 2.4**, using a common fixed-point method; we just implemented $\hat{\mathbf{m}}^{(\mathrm{Iter})} = F\left(\hat{\mathbf{m}}^{(\mathrm{Iter}-1)}\right)$ in (16) starting from an initial guess $\hat{\mathbf{m}}^{(0)}$ with the superscript being the number of $\mathrm{Iter} \in \mathbb{N}$, which always converges due to **Assumption 2** and (17). In this test case, we also consider the exponential moments $\mathbb{E}\left[e^{qX_t}\right]$ of $X$ and compare their values between the exact discretization methods and the generalized Riccati equation. For $q < 0$, the discretized Riccati equation is as follows: for $n = 0, 1, 2, \ldots$,

$$F_{(n+1)\Delta t'} = F_{n\Delta t'} + \frac{1}{N}\sum_{i=1}^{N} r_i \left\{ b_i G_{n\Delta t'}^{(i)} + a\int_0^{+\infty} \left(1 - e^{-G_{n\Delta t'}^{(i)} z_i}\right) v\left(\mathrm{d}z_i\right) \right\} \tag{69}$$

and

$$G_{(n+1)\Delta t'}^{(i)} = \frac{1}{1 + r_i \Delta t'} \left\{ G_{n\Delta t'}^{(i)} - r\frac{1}{2}\sigma^2 r_i \left(G_{n\Delta t'}^{(i)}\right)^2 + r_i \frac{1}{N}\sum_{j=1}^{N} A_{j,i} \int_0^{+\infty} \left(1 - e^{-G_{n\Delta t'}^{(j)} z_j}\right) v\left(\mathrm{d}z_j\right) \right\},\quad i = 1,2,3,\ldots,N, \tag{70}$$

where $\Delta t'$ is the time increment for discretizing the generalized Riccati equation. Therefore, the discretized equations (69) and (70) are based on the finite-dimensional version and classical forward Euler method for (69) and the semi-implicit Euler method for (70). The discretization for $q > 0$ is similar. The assumption about $q$ in **Proposition 1** is a sufficient condition but may not be necessary and sufficient. For each $q$, we empirically checked that choosing $\Delta t' = 0.001$ does not generate spurious oscillations of numerical solutions in numerical solutions to the generalized Riccati equations.

### 4.2.2 Results and discussion

Test 1 serves as the accuracy check of the discretization method in time and sampling duration because the average and variance do not depend on $\rho$ as demonstrated in (68). **Tables 5** and **6** summarize the computed average and variance for the moderate memory ($\alpha = 2.5$) and long memory cases ($\alpha = 0.8$), respectively. In both cases, the computed average and variance sampled in the longer duration are closer to the exact ones, suggesting that the sampling duration of 400,000 may not be necessarily sufficient for the long memory case. Indeed, it has been reported that excessively long sampling duration or burn-in period may become necessary for long memory processes [94,95]; our results are in accordance with these studies. To see the influences of the resolution in time $\Delta t$, **Table 7** lists the comparison of the computed and exact average and variance for $\alpha = 0.8$ and $N = 256$ with the sampling duration being the longer one. For both the average and variance, the relative error with respect to $\Delta t$ is nonmonotone and thus it appears to be difficult to estimate the convergence speed; however, it becomes larger for the coarser resolution.

**Figure 3** shows part of the computed sample paths of the superposition process $X$. **Figure 4** shows the corresponding autocorrelation functions based on sampling during the time interval $[1200000,\ 1600000]$. The fit between the computed and theoretical autocorrelation functions improves as $\alpha$ increases, that is, memory reduces. As listed in **Table 7**, the least-squares error between the computed and exact autocorrelation functions for the lag $0 \le \tau \le 100$ is one order of magnitude larger in this long-memory case ($\alpha \le 1$) than in the moderate-memory case ($\alpha > 1$). Moreover, the error is larger for a larger $N$ in some cases, which is considered to be owing to the existence of larger timescales with a larger $N$, the fluctuations of which persist more strongly in time. The agreement between the theoretical and computed autocorrelation functions shows that our methods work by properly reproducing the memory.

**Tables 5-8** do not imply clear convergence trends from the finite-dimensional model to the infinite-dimensional one as $N$ grows, although the size of errors seems not to be critical. This is considered due to at least the following three possible reasons. First, the computational error due to the time

discretization may have been comparable to the spatial one. According to **Table 7**, the error due to time discretization decreases as $N$ grows and is estimated to be the order of $O\left(10^{-3}\right)$ to $O\left(10^{-2}\right)$ when $\Delta t = O\left(10^{-2}\right)$. Second, errors due to the sampling method may also have affected the computational results. We estimated the statistics by using time averages, which would become more accurate as the memory becomes shorter where numerical solutions at consecutive time steps can be better approximated to be independent; however, this is not the case for our model because of its long-memory nature. Third, related to the second reason, growing $N$ implies not only finer discretization for $r > 0$ but also the introduction of smaller $r_i$ with which finite-dimensional nominal model involves slower timescales and hence slower numerical convergence of its statistics in time (see also **Figure 5**). Introducing smaller $r_i$ is not avoidable because smaller $r_i$ contributes to longer memory. This poses a severe computational cost in the Monte-Carlo simulation. **Tables 1-4** imply that the resolution in time and the total number of sample paths should be balanced even for the simpler square-root diffusion case. Derivation of the proper scaling relationship between the error of time discretization and that of sampling for the nominal model seems to be a nontrivial task for long-memory processes. Explicitly quantifying the errors due to finite $N$ would open a door towards the analysis in this direction.

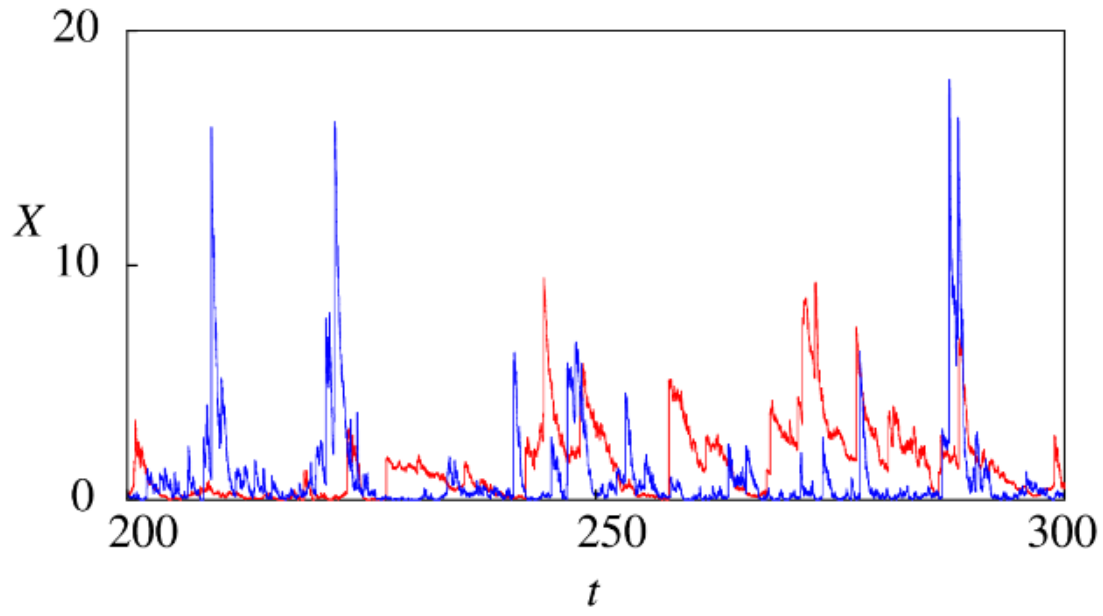

**Figure 3.** Computed sample paths of $X$ with $N = 256$. Colors: $\alpha = 0.8$ (red) and $\alpha = 2.5$ (blue).

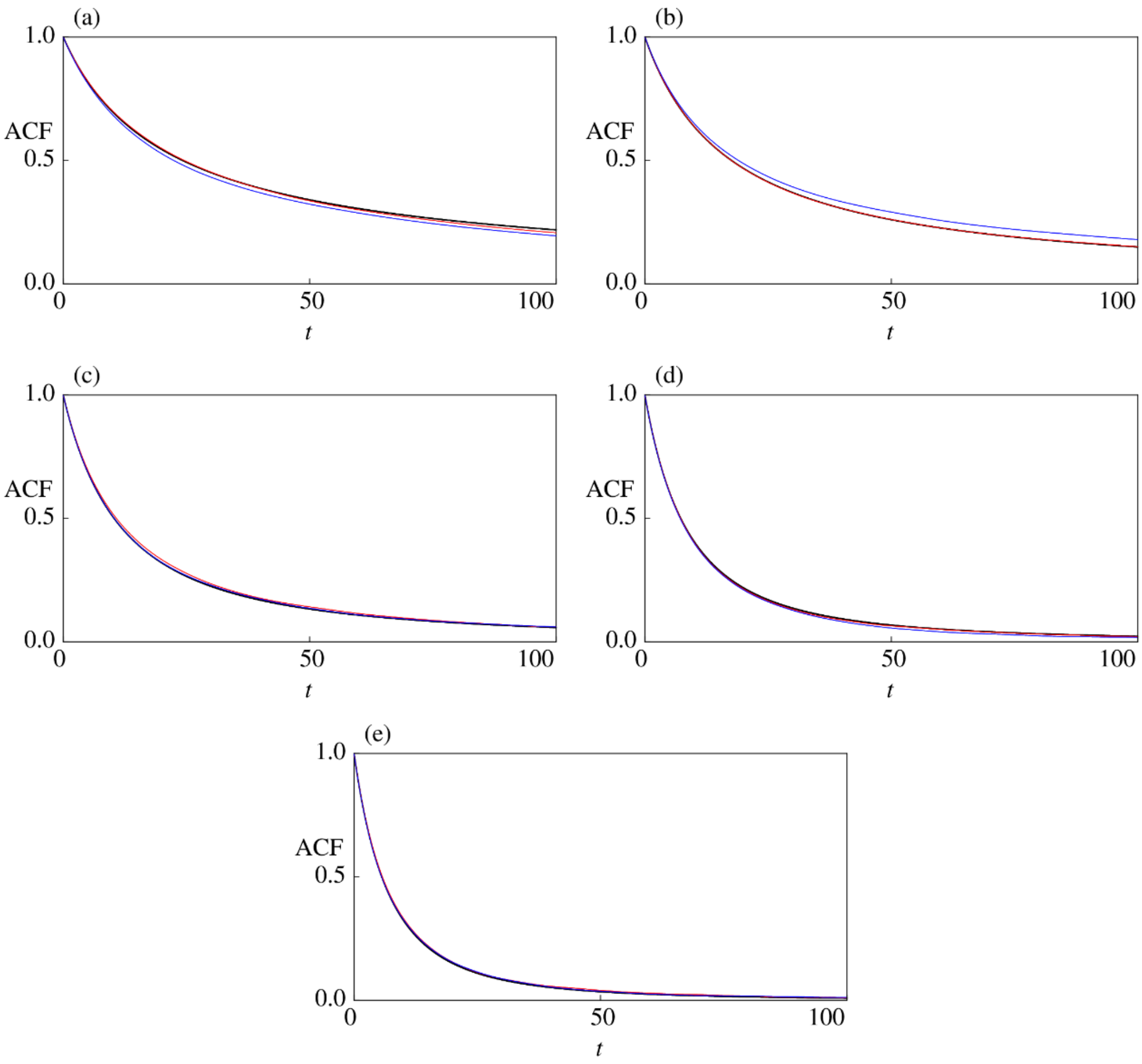

**Figure 4.** Comparison of computed and exact autocorrelation function (ACF) of $X$. Panels: (a) $\alpha = 0.8$, (b) $\alpha = 1.0$, (c) $\alpha = 1.5$, (d) $\alpha = 2.0$, (e) $\alpha = 2.5$. Colors: exact (black), $N = 256$ (red), $N = 512$ (blue).

**Table 5.** The computed average and its relative error (the numbers inside "( )"). The exact average is 1.

| | $\alpha = 0.8$ | | $\alpha = 2.5$ | |
|---|---|---|---|---|
| Sampling duration | Short | Long | Short | Long |
| $N = 64$ | 9.888.E-01 (1.120.E-02) | 9.981.E-01 (1.876.E-03) | 1.009.E+00 (9.060.E-03) | 1.006.E+00 (5.830.E-03) |
| $N = 128$ | 9.820.E-01 (1.796.E-02) | 9.979.E-01 (2.093.E-03) | 1.002.E+00 (2.070.E-03) | 1.001.E+00 (1.110.E-03) |
| $N = 256$ | 9.823.E-01 (1.773.E-02) | 9.979.E-01 (2.079.E-03) | 9.988.E-01 (1.226.E-03) | 9.969.R-01 (3.065.E-03) |

**Table 6.** The computed variance and its relative error (the numbers inside "( )"). The exact variance is 3.75.

| | $\alpha = 0.8$ | | $\alpha = 2.5$ | |
|---|---|---|---|---|
| Sampling duration | Short | Long | Short | Long |
| $N = 64$ | 3.670.E+00 (2.143.E-02) | 3.689.E+00 (1.615.E-02) | 3.813.E+00 (1.686.E-02) | 3.802.E+00 (1.394.E-02) |
| $N = 128$ | 3.477.E+00 (7.283.E-02) | 3.658.E+00 (2.458.E-02) | 3.743.E+00 (1.839.E-03) | 3.746.E+00 (1.168.E-03) |
| $N = 256$ | 3.510.E+00 (6.390.E-02) | 3.804.E+00 (1.429.E-02) | 3.697.E+00 (1.403.E-02) | 3.673.E+00 (2.050.E-02) |

**Table 7.** The computed variance and its relative error (the numbers inside "( )") for $\alpha = 0.8$ and $N = 256$. The sampling duration is the long one.

| $\Delta t$ | Average | Variance |
|---|---|---|
| 0.01 | 9.979.E-01 (2.079.E-03) | 3.804.E+00 (1.429.E-02) |
| 0.04 | 9.966.E-01 (3.416.E-03) | 3.734.E+00 (4.201.E-03) |
| 0.16 | 9.903.E-01 (9.663.E-03) | 3.590.E+00 (4.276.E-02) |
| 0.32 | 1.008.E+00 (7.670.E-03) | 3.602.E+00 (3.934.E-02) |

**Table 8.** The least-squares error between the computed and exact autocorrelation functions for the lag $0 \le \tau \le 100$.

| | $N = 256$ | $N = 512$ |
|---|---|---|
| $\alpha = 0.8$ | 1.496.E-02 | 2.876.E-02 |
| $\alpha = 1.0$ | 4.006.E-03 | 1.495.E-02 |
| $\alpha = 1.5$ | 5.042.E-03 | 3.069.E-03 |
| $\alpha = 2.0$ | 1.913.E-03 | 3.868.E-03 |
| $\alpha = 2.5$ | 2.302.E-03 | 2.008.E-03 |

We move to Test 2, where the focus is on the average that is computationally available with arbitrarily high accuracy for each $N$ owing to (16) and the exponential moment identified from the generalized Riccati equation. In contrast to Test 1, the average depends on $N$, i.e., the degree of freedom in the Markovian lifts, and we expect that by increasing $N$, the computational results would be closer to the exact results. Therefore, for each fixed $N$, a statistic obtained from (16) or the generalized Riccati equation is referred to as "Reference," and we compared the computed and reference results for different values of $N$. The sampling durations are the same as those in Test 1, and we fixed $\Delta t = 0.01$.

**Tables 9** and **10** list the comparison of the computed and reference average and exponential moment $\mathbb{E}\left[e^{qX_t}\right]$, respectively. For the exponential moment $\mathbb{E}\left[e^{qX_t}\right]$, we set $N = 512$ and examine the values of $q$ such that $q \le b_v = 0.1$ because the moment with $q > b_v$ does not exist. The accuracy with the exact discretization methods improves as the duration became longer, as for Test 1. All the exponential moments are reproduced for both large ($q = b_v = 0.1$) and small values ($q = -1$) of $q$ despite that their magnitudes are approximately four-order different. Further increasing the degree of freedom $N$ would yield results closer to the one with $N = +\infty$, but are not examined here because the computational costs may become prohibitive, particularly for the smaller $\alpha$.

The analysis here also elucidates the solution of the second generalized Riccati equation. **Figure 5** shows the time evolution of the numerical solution to the second generalized Riccati equation with $q = 0.1$, $\alpha = 0.8$, and $N = 512$. The numerical solution decreases for $r > 0$ and decreases exponentially toward zero at each discretized $r$ at large $t$. Interestingly, the exponent for a large $t$ is synchronized to the slowest one ($r_1 = 1/(2N)$) for all discretized $r$. This is because of the form of the second generalized Riccati equation, whose solution is fully coupled through a common integral term (the last term in (30) or (39)). The same tendency applies to the other cases. Thus, the computational results of the generalized Riccati equation implies that, at least computationally, the convergence toward the stationary state of the superposition process becomes slower as $N$ increases, i.e., as $r_1$ decreases. Because the generalized Riccati equation is an (deterministic) intego-differential equation, we expect some convergence of numerical solutions by using a (naive)finite difference method like that we employed here, but we did not have proof. A potential difficulty would be the unboundedness of the domain and coefficients with respect to $r > 0$, which may violate the global Lipschitz continuity (e.g., Lemma 3.5 in Cox et al. [96]): a critical condition for well-posing and hence approximating generalized Riccati equations. Moreover, numerical solutions obtained by the proposed method are nonnegative as shown in **Figure 5** and theoretically indicated in **Proposition 1**, while in general there will be no such guarantees if a discretization method based on some projection, such as a Galerkin or truncation method [e.g., 66], is used. We leave the numerical convergence of the generalized Riccati equation to be open.

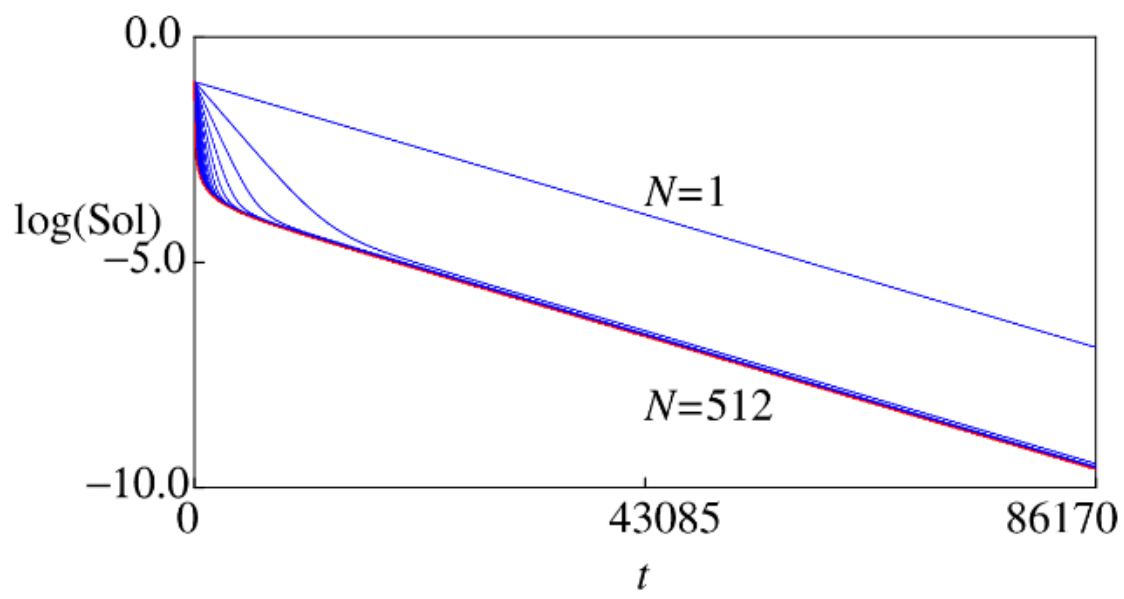

**Figure 5.** The numerical solution (Sol) to the second generalized Riccati equation at each discretized $r = r_i$ ($i = 1, 2, 3, ..., 512$). We plot the common logarithm of the solution (Sol) to visualize its exponential decay.

**Table 9.** Computed average of $X$ and its relative error (the numbers inside "( )"). "Reference" means the results due to solving the corresponding stationary linear system (16).

| | $\alpha = 0.8$ | | | $\alpha = 2.5$ | | |
|---|---|---|---|---|---|---|
| | Reference | Duration: Short | Duration: Long | Reference | Duration: Short | Duration: Long |
| $N = 64$ | 6.222.E+00 | 6.174.E+00 (7.671.E-03) | 6.214.E+00 (1.336.E-03) | 1.300.E+00 | 1.317.E+00 (1.349.E-02) | 1.308.E+00 (6.032.E-03) |
| $N = 128$ | 7.744.E+00 | 7.808.E+00 (8.256.E-03) | 7.751.E+00 (8.277.E-04) | 1.313.E+00 | 1.317.E+00 (3.374.E-03) | 1.316.E+00 (2.460.E-03) |
| $N = 256$ | 9.533.E+00 | 9.815.E+00 (2.966.E-02) | 9.736.E+00 (2.129.E-02) | 1.321.E+00 | 1.330.E+00 (7.338.E-03) | 1.325.E+00 (2.991.E-03) |

**Table 10.** Computed exponential moment $\mathbb{E}\left[e^{qX_t}\right]$ and its relative error (the numbers inside "( )"). "Reference" means the results due to the generalized Riccati equation.

| | $\alpha = 0.8$ | | $\alpha = 2.5$ | |
|---|---|---|---|---|
| | Reference | Computed | Reference | Computed |
| $q = 0.1$ | 4.119.E+00 | 4.059.E+00 (1.470.E-02) | 1.206.E+00 | 1.198.E+00 (6.659.E-03) |
| $q = 0.05$ | 1.868.E+00 | 1.865.E+00 (1.232.E-03) | 1.077.E+00 | 1.077.E+00 (4.178.E-04) |
| $q = -1$ | 6.914.E-04 | 7.128.E-04 (3.096.E-02) | 5.000.E-01 | 5.004.E-01 (7.180.E-04) |

### 4.3 Application to coupled model

The application case considers a coupled system of superposition processes that has recently been discussed by Yoshioka and Yoshioka [40] to show that exact discretization methods can deal with this advanced model. In this study, the model governs the river discharge $Q=(Q_t)_{t\in\mathbb{R}}$ and the non-seasonal part of the water quality index $C=(C_t)_{t\in\mathbb{R}}$ as follows: $Q$ is governed by the nominal model without diffusion ($b\equiv 0$) and self-exciting jumps ($a=1$, $A\equiv 0$) and $C$ by the nominal model without jumps ($a=0$, $A\equiv 0$) and discharge-dependent drift ($b=a_C+b_CQ_t$, $a_C,b_C>0$). Therefore, this model is a coupled version of two nominal models. We do not focus on the seasonality of water quality because it was modeled externally to the stochastic model.

Jump-driven discharge dynamics effectively reproduce flood events owing to runoff from the catchment of a river and discharge-driven dynamics associated with water quality changes [40, 97]. The discharge-dependent affine drift $a_C+b_CQ_t$ represents a positive correlation between water quality and quantity, which is often observed in indices such as total nitrogen and total organic carbon [98,99]. In this case, the model parameters are considered site-specific and identified using available time-series data. The fitted parameter values are summarized in Section 4 of Yoshioka and Yoshioka [40] ($a_C,b_C$ in this study is $a,\mu$ in the literature). These values were estimated using the moment-fitting method. In this study, sampling of the discharge $Q$ and concentration $C$ is conducted for 2,000 years after another 2,000 years to burn in. The time increment is fixed to 0.01 (day) unless otherwise specified, and we examine different values of $N$: 1024, 2048, and 4096.

**Table 11** lists the computed averages and variances of the discharge $Q$ and (non-seasonalized) total nitrogen concentration $C$ and their covariance, demonstrating that the results with the highest resolution $N$ provide the most accurate results in terms of the statistics, implying a certain convergence of the exact discretization methods. **Figure 6** shows the computed probability density functions (PDFs) of $Q$ and $C$, and **Figure 7** shows the joint PDFs. **Figure 6** shows that the computational results capture the tail of the discharge $Q$ but commonly fail to fit the peak of the mode in the PDF, the latter is owing to the fitted parameter values [40]. In contrast, concerning $C$, the peak of the empirical PDF is captured by all computational resolutions. Owing to the lack of empirical data, it is not possible to discuss the accuracy around the PDF tail of the computational results, but all predicted super-exponential decay. **Figure 7** suggests that the joint PDFs for different resolutions are comparable.

**Table 11.** Computed statistics of the coupled model: Average (Ave), Variance (Var), and Covariance (Cov).

| | Ave of $Q$ : ($m^3$/s) | Ave of $C$ : (-) | Var of $Q$ : ($m^6/s^2$) | Var of $C$ : (-) | Cov between $Q$ and $C$ : ($m^3$/s) |
|---|---|---|---|---|---|
| Theoretical | 4.156.E+01 | 1.084.E+00 | 2.755.E+03 | 2.676.E-01 | 9.778.E+00 |
| $N = 1024$ | 4.067.E+01 | 1.059.E+00 | 2.738.E+03 | 2.676.E-01 | 9.919.E+00 |
| | (2.127.E-02) | (2.302.E-02) | (5.971.E-03) | (1.607.E-04) | (1.438.E-02) |
| $N = 2048$ | 4.053.E+01 | 1.054.E+00 | 2.713.E+03 | 2.588.E-01 | 9.441.E+00 |
| | (2.464.E-02) | (2.751.E-02) | (1.529.E-02) | (3.294.E-02) | (3.448.E-02) |
| $N = 4098$ | 4.119.E+01 | 1.088.E+00 | 2.730.E+03 | 2.719.E-01 | 1.016.E+01 |
| | (8.896.E-03) | (3.967.E-03) | (9.176.E-03) | (1.613.E-02) | (3.912.E-02) |

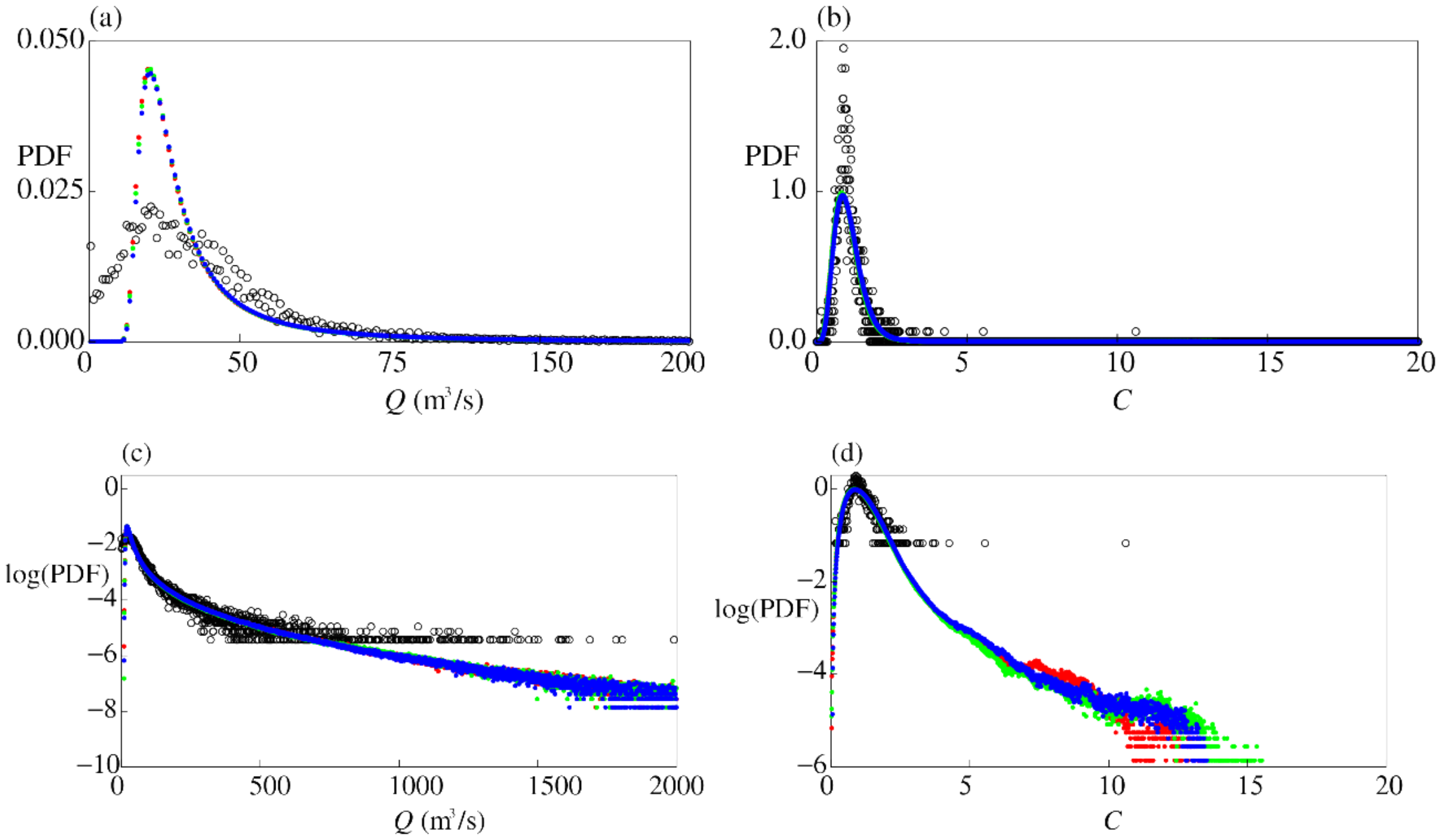

**Figure 6.** Computed PDFs in the logarithmic scale. Panels: (a) $Q$ in the ordinary scale, (b) $C$ in the ordinary scale, (c) $Q$ in the common logarithmic scale, and (d) $C$ in the common logarithmic scale. Symbols: Computed (blue) and theoretical (red). Colors: Empirical (black), $N = 1024$ (red), $N = 2048$ (green), $N = 4098$ (blue).

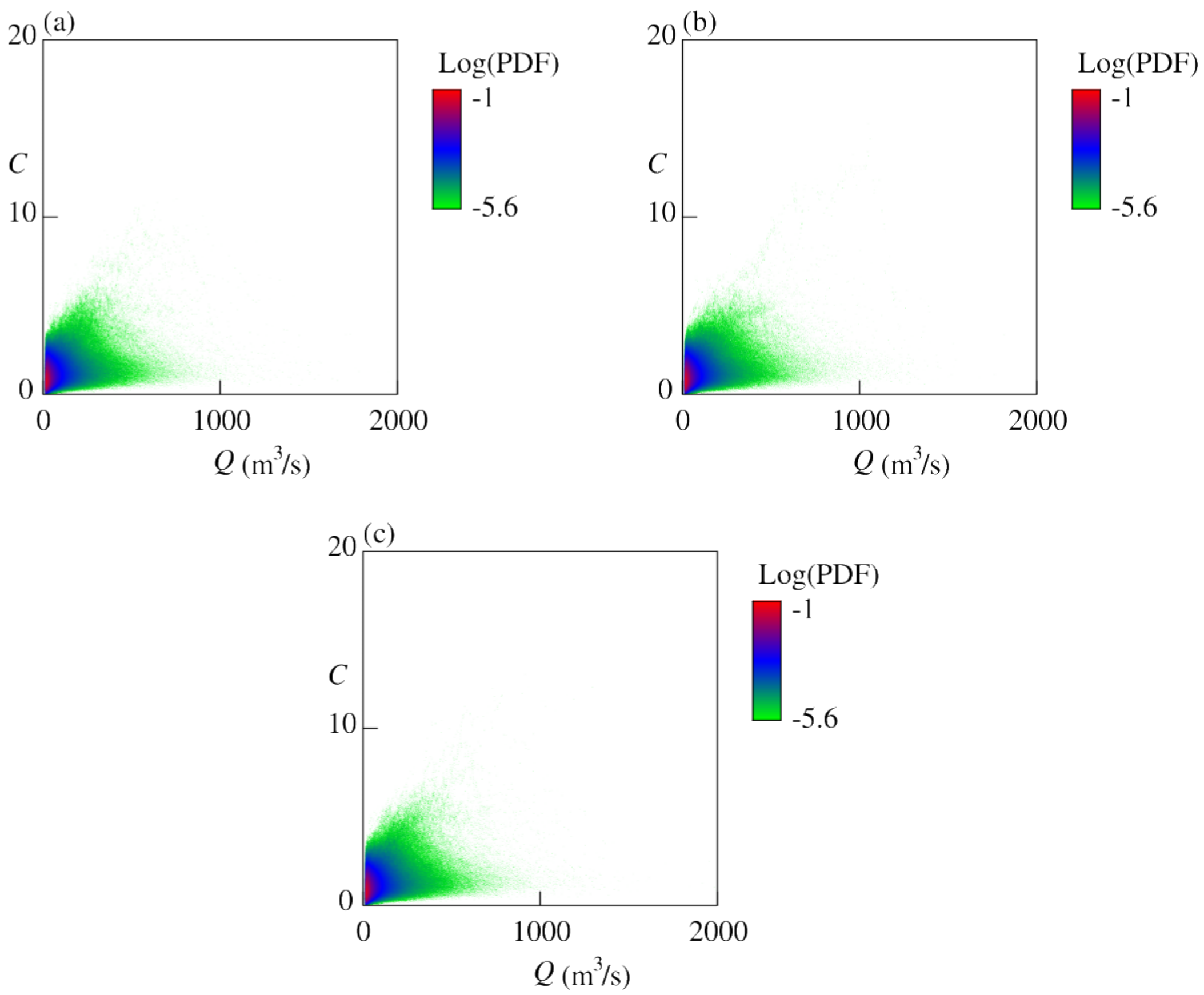

**Figure 7.** Computed joint PDFs of $Q$ and $C$ in the common logarithmic scale. Panels: (a) $N = 1024$, (b) $N = 2048$, (c) $N = 4096$.

### 4.4 Application to a nonlinear model

The last computational case shows that the proposed numerical method can handle non-affine superposition processes. We consider a CEV generalization of the superposition process $X$ in (1) with

$$\mathrm{d}x_t(S) = -rx_t(\mathrm{d}r)\mathrm{d}t + X_t^\gamma B(\mathrm{d}t,\mathrm{d}r) + \int_0^{+\infty} zJ(\mathrm{d}t,\mathrm{d}r,\mathrm{d}z) \quad (71)$$

with $\gamma \in \mathbb{R}$, which is reduced to the nominal model if $\gamma = 0$. The jump term is the same as that used in the second test case. The behavior of the classical CEV model without jumps has been summarized in previous studies based on the Feller condition [e.g., 100]; in the context of our model without superposition, the parameter $\gamma$ determines the hitting of solutions to the CEV model to boundary 0 such that there is no boundary hitting if $\gamma > 0$, always hit the boundary if $-1/2 < \gamma < 0$, and boundary hitting with large $\sigma$ if $\gamma = 0$. However, to the best of our knowledge, the superposition of CEV models has not yet been examined. We set $\sigma = 1$ and the values of the other parameters are the same as those in Test 2.

This generalization would allow for modeling a wider class of superposition processes that are genuinely nonlinear. Considering the corresponding finite-dimensional version as in **Section 2.3**, the diffusion coefficient of the process $x^{(i)}$ is set as

$$\left(X_t^{(N)}\right)^\gamma \sigma\sqrt{r_i x_t^{(i)}} = \sigma\sqrt{r_i}\frac{\sqrt{x_t^{(i)}}}{\left(\sum_{i=1}^N x_t^{(i)}\right)^{-\gamma}}, \quad (72)$$

which is bounded at the origin $x_t^{(1)} = x_t^{(2)} = \ldots = x_t^{(N)} = 0$ if $\gamma \geq -1/2$. The denominator of (72) is expected to converge to $X$ as $N \to +\infty$. The average of $X$ is the same as the second test case; therefore, we can compare the theoretical and numerical averages. We fix $\Delta t = 0.01$ and $N = 512$ with the sampling period of 1,600,000.

**Figure 8** shows the sample paths of the superposition process $X$ for different values of $\gamma$ and $\alpha$, indicating that the superposition processes still inherit the boundary hitting behaviors of the CEV model, which can be visually better inferred in **Figure 8(c)** with $\alpha = 2.5$, where the paths with $\gamma = -1/4$ and $\gamma = -1/2$ become closer to zero. **Figure 9** compares the computed PDFs of $X$ for different values of $\gamma$ and $\alpha$, suggesting longer tails for larger $\gamma$. **Table 12** lists the average and variance of $X$ for different values of $\gamma$ and $\alpha$, suggesting that the fluctuation of the average is larger among the cases with a smaller $\alpha = 0.8$ than those with a larger $\alpha = 2.5$, possibly because of the persistent temporal nonequilibrium fluctuation in the former, as in Test 1.

Finally, we analyze the influence of the time increment $\Delta t$. We fix $\alpha = 2.5$. **Table 13** lists the comparison of the computed average and variance of $X$ for different values of $\Delta t$ with $\gamma = 1/2$ and $\gamma = -1/2$. In both cases, the reference average is reproduced within a relative error of 1%, whereas the variance fluctuation is underestimated several times for a large $\Delta t$ when $\gamma = 1/2$. **Figure 10** shows the comparison of the computed PDFs of $X$ for different values of $\Delta t$, suggesting that the difference between the computed variances is owing to the extreme values constituting the tails of the computed PDFs, whereas their peaks are not critically different. Therefore, choosing a sufficiently small $\Delta t$ is necessary to reasonably estimate the tails of the PDFs in the nonlinear model.

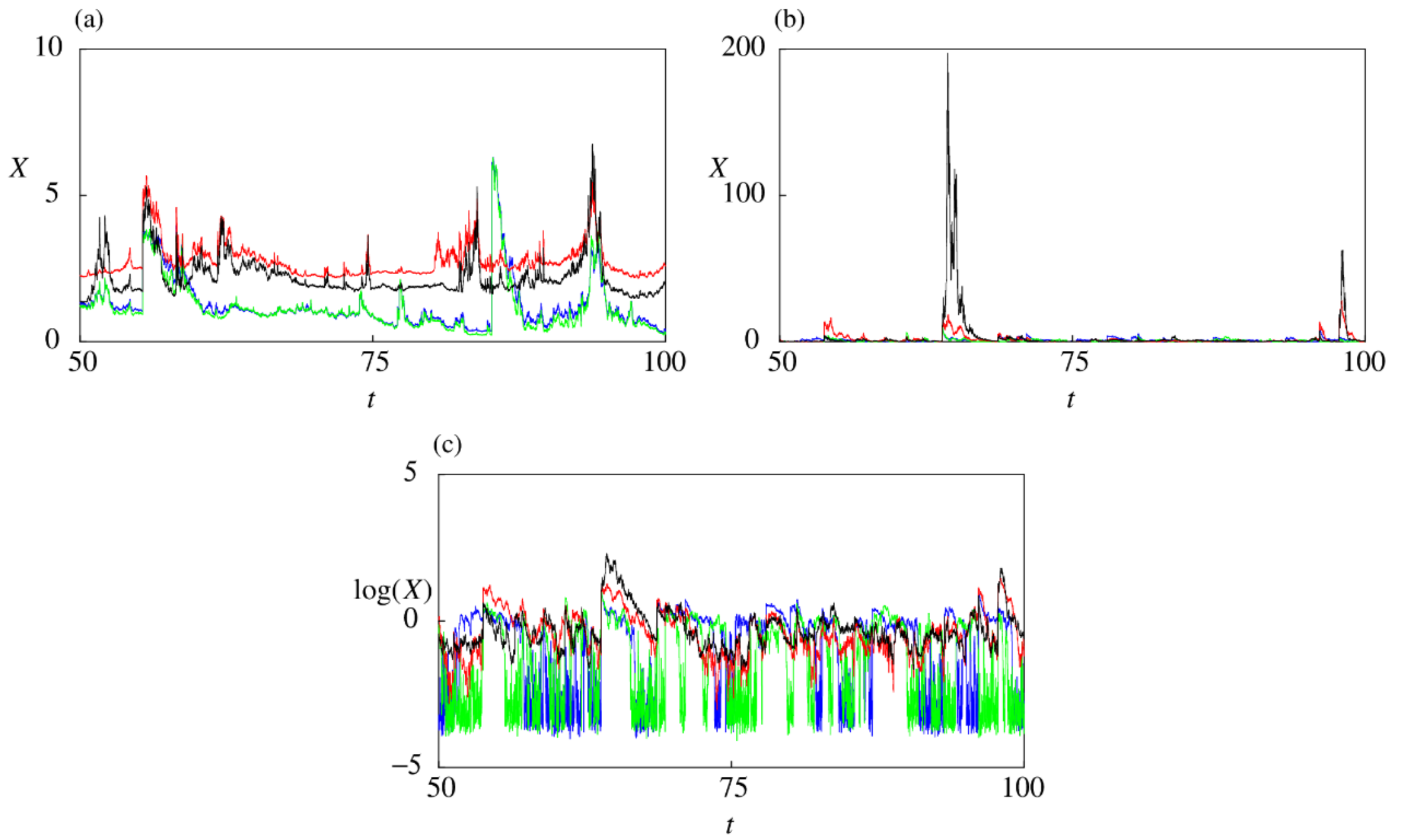

**Figure 8.** Computed sample paths of $X$. Panels: (a) $\alpha = 0.8$, (b) $\alpha = 2.5$, (c) $\alpha = 2.5$ in the common logarithmic scale. Colors: $\gamma = 1/2$ (black), $\gamma = 0$ (red), $\gamma = -1/4$ (green), $\gamma = -1/2$ (blue).

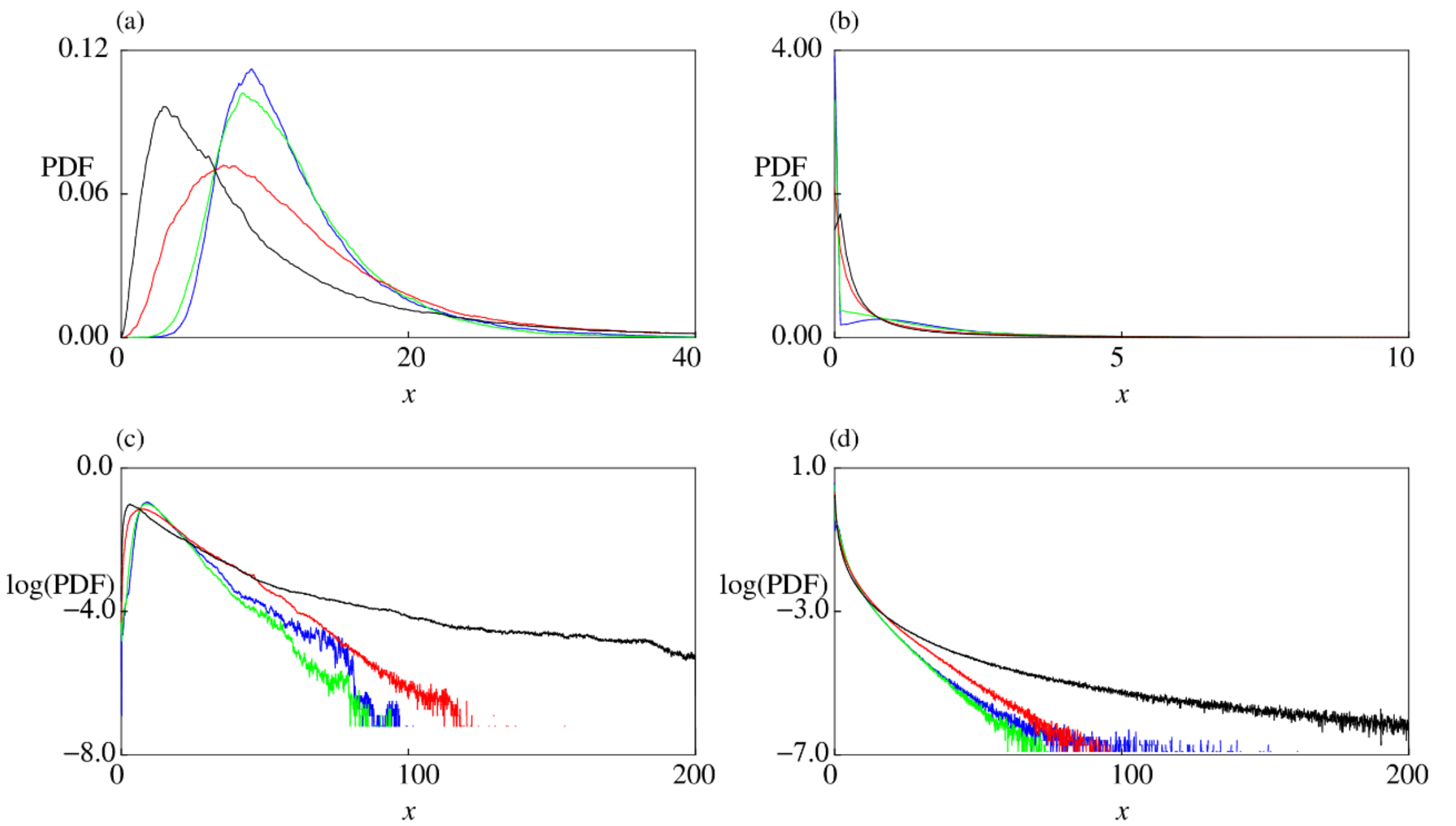

**Figure 9.** Computed PDFs of $X$. Panels: (a) $\alpha = 0.8$, (b) $\alpha = 2.5$, (a) $\alpha = 0.8$ in the common logarithmic scale, (d) $\alpha = 2.5$ in the common logarithmic scale. Colors: $\gamma = 1/2$ (black), $\gamma = 0$ (red), $\gamma = -1/4$ (green), $\gamma = -1/2$ (blue).

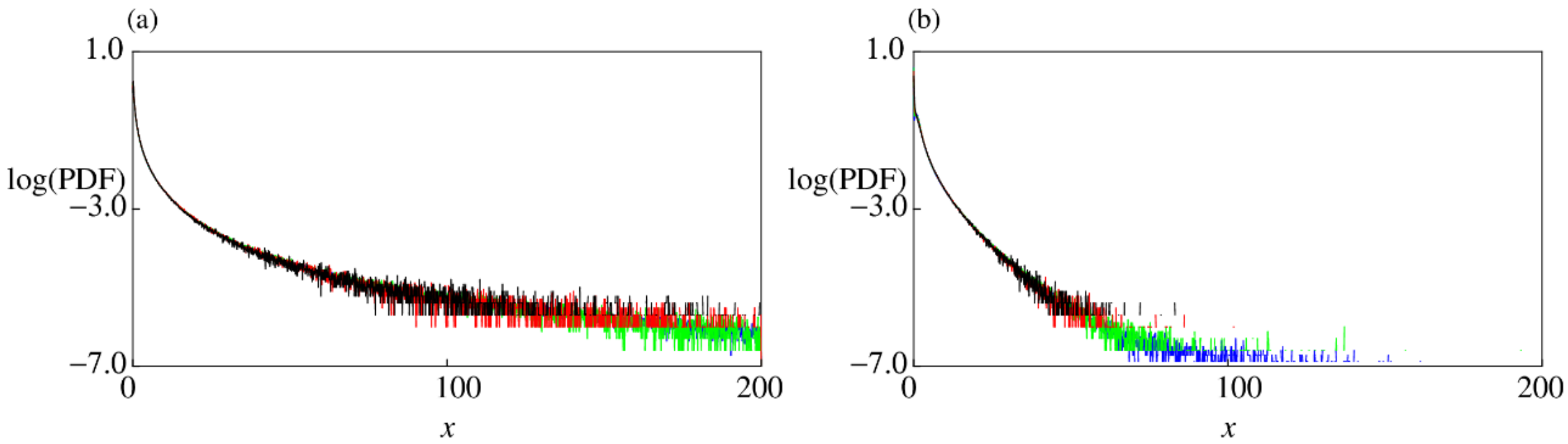

**Figure 10.** Computed PDFs of $X$ in the common logarithmic scale. Panels: (a) $\gamma = 1/2$, (b) $\gamma = -1/2$, (a) $\gamma = 1/2$. Colors: $\Delta = 0.01$ (black), $\Delta = 0.04$ (red), $\Delta = 0.16$ (green), $\Delta = 0.32$ (blue).

**Table 12.** Computed average and variance of $X$ for different values of $\gamma$. The numbers in "( )" represent relative errors.

| | $\alpha = 0.8$ | | $\alpha = 2.5$ | |
|---|---|---|---|---|
| | Average | Variance | Average | Variance |
| Reference | 11.6449 | | 1.32536 | |
| $\gamma = 1/2$ | 10.7818 (7.411.E-02) | 192.742 | 1.31757 (5.878.E-03) | 36.6101 |
| $\gamma = 0$ | 12.2359 (5.075.E-02) | 69.1488 | 1.32544 (6.036E-05) | 7.38221 |
| $\gamma = -1/4$ | 11.9397 (2.531.E-02) | 28.5423 | 1.32356 (1.358.E-03) | 5.27649 |
| $\gamma = -1/2$ | 12.2705 (5.372.E-02) | 33.2707 | 1.32212 (2.445.E-03) | 5.70047 |

**Table 13.** Computed average and variance of $X$ for different values of $\Delta t$. The numbers in "( )" represent relative errors.

| | $\gamma = 1/2$ | | $\gamma = -1/2$ | |
|---|---|---|---|---|
| | Average | Variance | Average | Variance |
| Reference | 1.32536 | | 1.32536 | |
| $\Delta t = 0.01$ | 1.31757 (5.878.E-03) | 36.6101 | 1.32212 (2.445.E-03) | 5.70047 |
| $\Delta t = 0.04$ | 1.32062 (3.576.E-03) | 26.4795 | 1.3261 (5.583.E-04) | 5.55538 |
| $\Delta t = 0.16$ | 1.33286 (5.659.E-03) | 27.4070 | 1.3249 (3.471.E-04) | 5.32841 |
| $\Delta t = 0.32$ | 1.33211 (5.093.E-03) | 24.6609 | 1.32317 (1.652.E-03) | 5.28742 |

## 5. Conclusions

We studied a nominal non-Markov model as a superposition process. The model was affine, driven by diffusion and jump noises, and accounted for the self-exciting nature to which a generalized Riccati equation was associated. The unique existence of solutions to the generalized Riccati equation depending on the initial conditions was then studied. We also addressed the computation of a nominal model based on a numerical scheme that exploited exact discretization methods for certain affine processes. The computational performance of the scheme was examined against test and application cases, demonstrating that it could handle superposition processes driven by jump and diffusion noise.

The outcomes of this study will open a new door to the modeling and computation of non-Markov processes of the superposition type, where the processes to be superposed may be nonlinear and neither linear nor affine. However, there are many theoretical issues concerning the well-posedness of nonlinear cases that should therefore be elaborated on in the future, including the global existence and regularity of solutions, although these issues can be resolved to some extent for finite-dimensional superpositions. One interesting issue to be investigated in the future is the study of a neural SDE version of the nominal and related models, where coefficients such as $b$ and $A$ in the nominal model are represented by a neural network [e.g., 101]. A more practical non-Markov process can be established if these coefficients are learned from time-series data. From a data science perspective, these coefficients can be learned from data. In such a case, sampling paths of the target process will become mandatory, and the tailored exact-discretization methods would be a candidate for computationally addressing this issue. From a computational perspective, the lack of a theoretical estimate of the convergence rate over time is a limitation of this study. We conjecture that the convergence rate may depend strongly on the parameter values, as in the classical square root case [e.g., 56]. Clarifying this point requires more effort from theoretical numerical analyses. Studying the convergence of numerical solutions to generalized Riccati equations will also be a key topic for better understanding superposition processes.

## Appendix: Proofs

***Proof of Lemma 1***

An elementary calculation shows that $\bar{\varphi}(q)$ is expressed as follows:

$$\bar{\varphi}(q)=\frac{q}{1-\bar{A}M_1-\frac{\sigma^2}{2}q}>0\,. \tag{73}$$

We can also check $\varphi<\gamma(\varphi)$ for $\varphi<\bar{\varphi}(q)$, and $\varphi>\gamma(\varphi)$ for $\varphi>\bar{\varphi}(q)$. Moreover, (73) shows $\bar{\varphi}(q)>q$ and $\lim_{q\to+0}\bar{\varphi}(q)=0$.

□

***Proof of Proposition 1***

We split the proof into three steps as follows: The strategy of the proof is qualitatively based on that of Yoshioka [102]; however, several modifications are necessary because our model is more complex.

At the first step, we study the auxiliary equation

$$\varphi(t,r)=\mathbb{G}(\hat{\varphi})(t,r)\left(:=\hat{\mathbb{G}}(\varphi)(t,r)\right),\ t\geq 0\ \text{and}\ r>0 \tag{74}$$

with the truncation function $\hat{\varphi}(\cdot,\cdot)=\max\left\{0,\min\left\{\bar{\varphi}(q),\varphi(\cdot,\cdot)\right\}\right\}$ (the upper truncation $\bar{\varphi}(q)$ comes from **Lemma 1**), which follows by the boundedness, Lipschitz continuity, and strict contraction properties of $\hat{\mathbb{G}}$. For each $T>0$, it holds true that $\hat{\mathbb{G}}$ is a mapping from $C_{\mathrm{b},T}$ to $C_{\mathrm{b},T}$, because $qe^{-rt}$ is bounded and continuous , $A$ is bounded and continuous , and the following integrals for $t\in[0,T]$ and $r>0$ are continuous and bounded: $r\int_{s=0}^{s=t}e^{-\eta(q)r(t-s)}\mathrm{d}s$, $r\int_{s=0}^{s=t}e^{-\eta(q)r(t-s)}\left(\bar{\varphi}(q)-\hat{\varphi}(t,r)\right)\hat{\varphi}(t,r)\mathrm{d}s$, and $r\int_{s=0}^{s=t}e^{-\eta(q)r(t-s)}\int_{u=0}^{u=+\infty}A(u,r)\int_{z=0}^{z=+\infty}\left(1-e^{-\hat{\varphi}(s,u)z}\right)v(\mathrm{d}z)\rho(\mathrm{d}u)\mathrm{d}s$. For boundedness, also refer to (75) and (76). Continuity follows from the dominated convergence theorem (Theorem 4.20 in Farenick [103]) because the integrands of these integrals are compositions of continuous functions that can be bounded by integrable functions, that is, $0\leq e^{-r(t-s)}\leq 1$, $0\leq\left(\bar{\varphi}(q)-\hat{\varphi}(t,r)\right)\hat{\varphi}(t,r)\leq\left(\bar{\varphi}(q)\right)^2$, $0\leq A(u,r)\leq\bar{A}$, and $0\leq 1-e^{-\hat{\varphi}(s,u)z}\leq\hat{\varphi}(s,u)z\leq\bar{\varphi}(q)z$. The range of integration $(0,t)$ in these integrals can be transformed into $(0,1)$ by introducing the new auxiliary variable $s'=ts$ without losing the boundedness and continuity of their integrands.

We first show that the auxiliary equation (74) admits a unique solution, denoted by $\tilde{\varphi}$. Second, we show that $\tilde{\varphi}$ satisfies the global bound $0\leq\tilde{\varphi}\leq\bar{\varphi}(q)$ stated in the proposition, indicating that a solution exists for the second Riccati equation (30) that satisfies the required bound. Finally, the strict contraction property of $\hat{\mathbb{G}}$ in $C_{\mathrm{b},T}$ shows that the unique solution to the auxiliary equation (74) is the desired solution to the second Riccati equation (30). We use the following elementary inequalities: $1-e^{-x} \le x$ and $\left|e^{-x}-e^{-y}\right| \le |x-y|$ for $x, y \ge 0$, and

$$\left\|\int_{s=0}^{s=t} e^{-\eta(q)r(t-s)} r\mathrm{d}s\right\|_T = \sup_{0\le t\le T}\sup_{r>0}\int_{s=0}^{s=t} e^{-\eta(q)r(t-s)} r\mathrm{d}s \le \frac{1}{\eta(q)}\sup_{0\le t\le T}\sup_{r>0}\left(1-e^{-rt}\right) \le \frac{1}{\eta(q)} \quad \text{for any } T>0. \quad (75)$$

***First step***

The boundedness and Lipschitz continuity of $\hat{\mathbb{G}}$ are proved as follows. Fix arbitrary $T > 0$, and choose arbitrary $\varphi_1, \varphi_2 \in C_{\mathrm{b},T}$. We have, by the triangle inequality,

$$\begin{aligned}
\left\|\hat{\mathbb{G}}(\varphi_1)\right\|_T &\le q + \frac{1}{2}\sigma^2 \left\|\int_{s=0}^{s=t} e^{-\eta(q)r(t-s)} r\left(\bar{\varphi}(q)-\hat{\varphi}_1(s,r)\right)\hat{\varphi}_1(s,r)\mathrm{d}s\right\|_T \\
&\quad + \left\|\int_{s=0}^{s=t} e^{-\eta(q)r(t-s)} r\int_{u=0}^{u=+\infty} A(u,r)\int_{z=0}^{z=+\infty}\left(1-e^{-\hat{\varphi}_1(s,u)z}\right)v(\mathrm{d}z)\rho(\mathrm{d}u)\mathrm{d}s\right\|_T \\
&\le q + \frac{1}{2}\sigma^2\left(\bar{\varphi}(q)\right)^2 \left\|\int_{s=0}^{s=t} e^{-\eta(q)r(t-s)} r\mathrm{d}s\right\|_T \\
&\quad + \bar{A}\left\|\int_{s=0}^{s=t} e^{-\eta(q)r(t-s)} r\int_{u=0}^{u=+\infty}\int_{z=0}^{z=+\infty}\left(1-e^{-\bar{\varphi}(q)z}\right)v(\mathrm{d}z)\rho(\mathrm{d}u)\mathrm{d}s\right\|_T \\
&\le q + \frac{1}{2}\sigma^2\left(\bar{\varphi}(q)\right)^2 \left\|\int_{s=0}^{s=t} e^{-\eta(q)r(t-s)} r\mathrm{d}s\right\|_T + \bar{A}\left\|\int_{s=0}^{s=t} e^{-\eta(q)r(t-s)} r\int_{z=0}^{z=+\infty}\left(\bar{\varphi}(q)z\right)v(\mathrm{d}z)\mathrm{d}s\right\|_T \qquad (76) \\
&= q + \frac{1}{2}\sigma^2\left(\bar{\varphi}(q)\right)^2 \left\|\int_{s=0}^{s=t} e^{-\eta(q)r(t-s)} r\mathrm{d}s\right\|_T + \bar{A}M_1\bar{\varphi}(q)\left\|\int_{s=0}^{s=t} e^{-\eta(q)r(t-s)} r\mathrm{d}s\right\|_T \\
&\le q + \left(\frac{1}{2}\sigma^2\left(\bar{\varphi}(q)\right)^2 + \bar{A}M_1\bar{\varphi}(q)\right)\frac{1}{\eta(q)} \\
&= q + \bar{\varphi}(q)\left(\bar{A}M_1 + \frac{1}{2}\sigma^2\bar{\varphi}(q)\right)\left(1+\frac{1}{2}\sigma^2\bar{\varphi}(q)\right)^{-1} \\
&= \bar{\varphi}(q),
\end{aligned}$$

where we used **Lemma 1** to obtain the last line. For Lipschitz continuity, we have, again by the triangle inequality,

$$
\begin{aligned}
&\left\|\hat{\mathbb{G}}(\varphi_1)-\hat{\mathbb{G}}(\varphi_2)\right\|_T \\
&\leq \frac{1}{2}\sigma^2\bar{\varphi}(q)\left\|\int_{s=0}^{s=t} e^{-\eta(q)r(t-s)} r\left|\hat{\varphi}_1(s,r)-\hat{\varphi}_2(s,r)\right|\mathrm{d}s\right\|_T \\
&+\bar{A}\left\|\int_{s=0}^{s=t} e^{-\eta(q)r(t-s)} r\int_{u=0}^{u=+\infty}\int_{z=0}^{z=+\infty}\left|\hat{\varphi}_1(s,u)-\hat{\varphi}_2(s,u)\right| z v(\mathrm{d}z)\rho(\mathrm{d}u)\mathrm{d}s\right\|_T \\
&=\frac{1}{2}\sigma^2\bar{\varphi}(q)\left\|\int_{s=0}^{s=t} e^{-\eta(q)r(t-s)} r\left|\hat{\varphi}_1(s,r)-\hat{\varphi}_2(s,r)\right|\mathrm{d}s\right\|_T \\
&+\bar{A}M_1\left\|\int_{s=0}^{s=t} e^{-\eta(q)r(t-s)} r\left(\int_{u=0}^{u=+\infty}\left|\hat{\varphi}_1(s,u)-\hat{\varphi}_2(s,u)\right|\rho(\mathrm{d}u)\right)\mathrm{d}s\right\|_T \\
&\leq \frac{1}{2}\sigma^2\bar{\varphi}(q)\left\|\int_{s=0}^{s=t} e^{-\eta(q)r(t-s)} r\left|\hat{\varphi}_1(s,r)-\hat{\varphi}_2(s,r)\right|\mathrm{d}s\right\|_T \\
&+\bar{A}M_1\left\|\int_{s=0}^{s=t} e^{-\eta(q)r(t-s)} r\left(\int_{u=0}^{u=+\infty}\left\|\hat{\varphi}_1(s,\cdot)-\hat{\varphi}_2(s,\cdot)\right\|\rho(\mathrm{d}u)\right)\mathrm{d}s\right\|_T \\
&\leq \frac{1}{2}\sigma^2\bar{\varphi}(q)\left\|\int_{s=0}^{s=t} e^{-\eta(q)r(t-s)} r\left\|\varphi_1(s,\cdot)-\varphi_2(s,\cdot)\right\|\mathrm{d}s\right\|_T \\
&+\bar{A}M_1\left\|\int_{s=0}^{s=t} e^{-\eta(q)r(t-s)} r\left\|\varphi_1(s,\cdot)-\varphi_2(s,\cdot)\right\|\mathrm{d}s\right\|_T \\
&\leq\left(\frac{1}{2}\sigma^2\bar{\varphi}(q)+\bar{A}M_1\right)\left\|\int_{s=0}^{s=t} e^{-\eta(q)r(t-s)} r\left\|\varphi_1(s,\cdot)-\varphi_2(s,\cdot)\right\|\mathrm{d}s\right\|_T .
\end{aligned}
\tag{77}
$$

We also have

$$
\begin{aligned}
\left\|\int_{s=0}^{s=t} e^{-\eta(q)r(t-s)} r\left\|\varphi_1(s,\cdot)-\varphi_2(s,\cdot)\right\|\mathrm{d}s\right\|_T &= \sup_{0\leq t\leq T}\left\|\int_{s=0}^{s=t} e^{-\eta(q)r(t-s)} r\left\|\varphi_1(s,\cdot)-\varphi_2(s,\cdot)\right\|\mathrm{d}s\right\|_T \\
&\leq \sup_{0\leq t\leq T}\left(\left\|\int_{s=0}^{s=t} e^{-\eta(q)r(t-s)} r\left\|\varphi_1-\varphi_2\right\|_T \mathrm{d}s\right\|\right) \\
&=\left\|\varphi_1-\varphi_2\right\|_T \sup_{0\leq t\leq T}\left(\left\|\int_{s=0}^{s=t} e^{-\eta(q)r(t-s)} r\mathrm{d}s\right\|\right) \\
&\leq \frac{1}{\eta(q)}\left\|\varphi_1-\varphi_2\right\|_T .
\end{aligned}
\tag{78}
$$

Consequently, we obtain the Lipschitz continuity as follows:

$$
\left\|\hat{\mathbb{G}}(\varphi_1)-\hat{\mathbb{G}}(\varphi_2)\right\|_T \leq\left(\frac{1}{2}\sigma^2\bar{\varphi}(q)+\bar{A}M_1\right)\frac{1}{\eta(q)}\left\|\varphi_1-\varphi_2\right\|_T . \tag{79}
$$

The coefficient multiplied by $\left\|\varphi_1-\varphi_2\right\|_T$ in (79) is smaller than 1 because of $\bar{A}M_1<1$, leading to the strict contraction property:

$$
\left(\frac{1}{2}\sigma^2\bar{\varphi}(q)+\bar{A}M_1\right)\frac{1}{\eta(q)}=\left(\frac{1}{2}\sigma^2\bar{\varphi}(q)+\bar{A}M_1\right)\left(\frac{1}{2}\sigma^2\bar{\varphi}(q)+1\right)^{-1}<1 . \tag{80}
$$

The mapping $\hat{\mathbb{G}}: C_{\mathrm{b},T}\to C_{\mathrm{b},T}$ is therefore strictly contractive, proving that by Banach fixed point theorem (e.g., Theorem (1.1) in Granas and Dugundji [104]) there is a unique solution to the auxiliary equation (74) in each time interval $[0,T]$ with arbitrary $T>0$. This solution is denoted by $\tilde{\varphi}$, which can be extended to an arbitrary $t>0$ because of the estimate (as in (76)), which is uniform in time:

$$
\left\|\varphi(t,\cdot)\right\|=\left\|\hat{\mathbb{G}}(\varphi)(t,\cdot)\right\|\leq\bar{\varphi}(q),\quad t>0 . \tag{81}
$$

***Second step***

We show $0 \le \tilde{\varphi}(t,r) \le \bar{\varphi}(q)$ for all $t \ge 0$ and $r > 0$. The upper bound of $\tilde{\varphi}$ is due to (81). The lower bound follows from

$$\begin{aligned}\tilde{\varphi}(t,r) &= \hat{\mathbb{G}}(\tilde{\varphi})(t,r)\\ &= qe^{-\eta(q)rt} + \int_{s=0}^{s=t} e^{-\eta(q)r(t-s)} r \frac{1}{2}\sigma^2 \left(\bar{\varphi}(q) - \tilde{\varphi}(s,r)\right)\tilde{\varphi}(s,r)\mathrm{d}s\\ &\quad + \int_{s=0}^{s=t} e^{-\eta(q)r(t-s)} r \int_{u=0}^{u=+\infty} A(u,r)\int_{z=0}^{z=+\infty}\left(1-e^{-\tilde{\varphi}(s,u)z}\right)v(\mathrm{d}z)\rho(\mathrm{d}u)\mathrm{d}s\\ &\ge 0\end{aligned} \tag{82}$$

because $q > 0$ and the integrand of (82) are nonnegative. Therefore, $\tilde{\varphi}$ solves (30).

***Last step***

We show that $\tilde{\varphi}$ is the unique nonnegative solution to (30) such that $\|\tilde{\varphi}\|_T \le \bar{\varphi}(q)$. Assume that there exist two such solutions $\varphi_1, \varphi_2 \in C_{\mathrm{b},T}$ to (30). Subsequently, we have (as for (77)–(79))

$$\|\varphi_1 - \varphi_2\|_T = \|\mathbb{G}(\varphi_1) - \mathbb{G}(\varphi_2)\|_T \le \left(\frac{1}{2}\sigma^2\bar{\varphi}(q) + \bar{A}M_1\right)\frac{1}{\eta(q)}\|\varphi_1 - \varphi_2\|_T < \|\varphi_1 - \varphi_2\|_T, \tag{83}$$

and thus, $\|\varphi_1 - \varphi_2\|_T = 0$, showing the uniqueness.

Finally, $\mathbb{G}(\tilde{\varphi})$ as a function of $t$ (for the solution $\tilde{\varphi}$ and each fixed $r > 0$) is continuously differentiable with respect to $t \in [0,T]$ by Proposition 2' in Zorich [105], because the integrand of $\mathbb{G}(\tilde{\varphi})$ is a continuous function of $(s,t) \in [0,T]$. This implies that the solution $\tilde{\varphi}$ is also continuously differentiable with respect to $t \in (0,T]$. The proof is then completed because $T > 0$ is arbitrary.

□

***Proof of Proposition 2***

For a sufficiently small $q > 0$, the integral on the right-hand side of (37) is uniformly bounded for $t > 0$ because $G$ is bounded from above by the constant $\bar{\varphi}(q)$. The integrand in time integration is strictly bounded:

$$\int_{r=0}^{r=+\infty} r\left\{b(r)G(s,r) + a\int_{z=0}^{z=+\infty}\left(1-e^{-G(s,r)z}\right)v(\mathrm{d}z)\right\}\rho(\mathrm{d}r) \le \left(\sup_{r>0} b(r) + aM_1\right)\bar{\varphi}(q)\int_{r=0}^{r=+\infty} r\rho(\mathrm{d}r) < +\infty \tag{84}$$

for any $s > 0$. Moreover, this integrand is continuous, which completes the proof.

□